\documentclass{amsart}
\usepackage{graphicx}
\usepackage{amssymb}
\usepackage{amsfonts}
\usepackage{hyperref}
\usepackage{mathrsfs}
\swapnumbers
\newtheorem{theorem}{Theorem}[section]
\newtheorem{lemma}[theorem]{Lemma}

\newtheorem{proposition}[theorem]{Proposition}
\theoremstyle{definition}

\newtheorem{assumption}[theorem]{Assumption}
\newtheorem{remark}[theorem]{Remark}

\numberwithin{equation}{section}
\theoremstyle{plain}

\numberwithin{equation}{section} 
\numberwithin{figure}{section} 
\theoremstyle{plain}
\theoremstyle{plain}
\theoremstyle{remark}
\newtheorem*{acknowledgement*}{Acknowledgement}
\theoremstyle{example}

\newcommand{\cC}{{\mathcal C}}

\newcommand{\cE}{{\mathcal E}}
\newcommand{\cF}{{\mathcal F}}

\newcommand{\cH}{{\mathcal H}}

\newcommand{\cK}{{\mathcal K}}
\newcommand{\cL}{{\mathcal L}}

\newcommand{\cP}{{\mathcal P}}

\newcommand{\cR}{{\mathcal R}}
\newcommand{\cS}{{\mathcal S}}

\newcommand{\te}{{\theta}}

\newcommand{\om}{{\omega}}
\newcommand{\ve}{{\varepsilon}}
\newcommand{\del}{{\delta}}
\newcommand{\Del}{{\Delta}}
\newcommand{\gam}{{\gamma}}
\newcommand{\Gam}{{\Gamma}}

\newcommand{\sig}{{\sigma}}
\newcommand{\al}{{\alpha}}
\newcommand{\be}{{\beta}}
\newcommand{\ka}{{\kappa}}
\newcommand{\la}{{\lambda}}
\newcommand{\Lam}{{\Lambda}}

\newcommand{\bbC}{{\mathbb C}}

\newcommand{\bbN}{{\mathbb N}}

\newcommand{\bbR}{{\mathbb R}}

\newcommand{\bbZ}{{\mathbb Z}}
\newcommand{\bbI}{{\mathbb I}}
\newcommand{\bbQ}{{\mathbb Q}}

\begin{document}
\title[]{Complex projective metrics on Young towers, countable sub-shifts and other countable covering  maps} 
 \vskip 0.1cm
 \author{Yeor Hafouta \\
\vskip 0.1cm
Einstein Institute  of Mathematics\\
Hebrew University\\
Jerusalem, Israel}%
\address{
Einstein Institute of Mathematics, The Hebrew University, Jerusalem 91904, Israel}
\email{yeor.hafouta@mail.huji.ac.il}%

\thanks{ }
\dedicatory{  }
 \date{\today}

\maketitle
\markboth{Y. Hafouta}{Complex cones for towers, shifts and covering maps} 
\renewcommand{\theequation}{\arabic{section}.\arabic{equation}}
\pagenumbering{arabic}

\begin{abstract}\noindent
In this article we  apply complex projective metrics to sequences of complex transfer operators generated by Young towers, countable shifts and other types of distance expanding maps (possibly time dependent) with countable degrees.
We will derive a sequential Ruelle-Perron-Frobenius theorem for such sequences, and relying on it a variety of probabilistic limit theorems (finer than the CLT) for non-stationary processes follow.
\end{abstract}

\section{Introduction}\label{sec1}
There are several techniques to prove that complex perturbations of real quasi compact operator are quasi-compact.
One of the most common techniques relies on application of perturbation theorems (see  \cite{HH}, \cite{Kato}). This method can be applied with a some classes of transfer operators,
(see, for instance \cite{HH}, \cite{GH},\cite{Jon} and \cite{Y1}), where the quasi compactness is  verified via a Lasota-Yorke type inequality, or by finite dimensional range approximation arguments (see \cite{Y1}). 
 This classical scheme does not work when dealing with compositions of different complex transfer operators (even if they are different perturbations of the same operator), since there is no spectral theory to exploit. 
 
Another approach for proving that a real (transfer) operator is quasi-compact relies on contraction properties of projective Hilbert metrics (see \cite{Bir}, \cite{Bus} and \cite{Liverani}). This method does not rely on spectral-theory arguments, and so it also works well (for certain models) in the case of composition of different real operators. For instance, it yields
 certain type of Ruelle-Perron-Frobenuius (RPF) theorem for some classes of (compositions of) stationary random transfer operators (taken along paths of measure preserving systems), which is a specific case of time dependent operators (see \cite{Kifer-1996}, \cite{Kifer-Thermo} and references therein). Such RPF theorems have a variety of applications, for instance to limit theorem for random dynamical system (see \cite{Kifer-1996}, \cite{Kifer-1998} and reference therein). Several types of limit theorems were also proved for time dependent sequential dynamical systems  using different methods (see \cite{Aimino}, \cite{Arno}, \cite{Conze}, \cite{Nicol} and references therein). All the above results  do not include fine limit theorems such as the Berry-Esseen theorem and the local central limit theorem (LCLT).

In \cite{Rug} (see also \cite{Dub1} and \cite{Dub2}) the theory of complex cones and complex projective Hilbert metrics was developed, and the contraction properties of linear operators between two complex cones were established. Moreover, a certain type of (conic) perturbation theory was developed.
In \cite{book} we used this theory and obtained a version of the RPF theorem for compositions of random complex transfer operators generated by a statioanry seqeunce of different distance exanding maps $T_j=T_{\te^j\om}$, and in \cite{SeqRPF} we extended it to the case of a single sequence $T_j$ (i.e. which  is not necessarily random), and also studied some type of stability properties of the resulting RPF triplets. These RPF theorems have many applications:
we applied them in order to derive several limit theorems (e.g. LCLT and Berry-Esseen theorems) for large classes of non-stationary processes (see \cite{book}, \cite{HafEdge} and \cite{SeqRPF}) generated by a sequnece of maps $T_j$, for certain classes of skew products (\cite{Annealed}) and for so called, nonconventional sums (see Ch. 2 of  \cite{book}). All of these results were applicable only for transfer operators generated by finite degree sequence of maps $T_j$.

In this paper we  apply complex cone techniques for countable degree maps, and obtain sequential complex RPF theorems for several types of sequences of complex perturbations  of an underlying sequence of real transfer operators $\cL_0^{(j)}$. 
We first consider a sequence $\cL_z^{(j)}$ of complex perturbations of a singe transfer operator $\cL_0=\cL_0^{(j)}$ defined by a Young tower or a countable subshift of finite type. More specifically, we  consider the deterministic quasi-compact transfer operators $\cL_0$ from \cite{Viv1} and \cite{Viv2} and then consider perturbations of the form $\cL_z^{(j)}=\cL_0(e^{zu_j})$, where $u_j$ is a sequence of H\"older continuous functions. Exactly as in \cite{SeqRPF}, this complex sequential RPF theorem yields various limit theorems (finer than the central limit theorem) for random variables of the form $\sum_{k=0}^{n-1}u_j(T^jx)$, where $T$ is either a Young tower or a countable shift, and $x$ is distributed according to an appropriate conformal or invariant measure. 

We also return to the setup of sequential covering maps $T_j$, were, in contrast to \cite{SeqRPF}, we will allow  that some (possibly all) of the $T_j$'s will have countable degree. The RPF theorem we prove in this setup holds true for locally distance expanding interval maps, including the ones considered in Section 3 of \cite{Conze}, and also in the case when each $T_j=T$ is the Gauss map and $\cL_0$ is the classical transfer operator given by 
\[
\cL_0g(x)=\sum_{n=1}^\infty\frac1{(n+x)^2}g\big(\frac1{(n+x)^2}\big).
\]
Our results also hold true for countable nonstationary full-shifts, and also when each $T_j$ is the same countable subshift of finite type which has a Gibbs measure (see \cite{Sarig1}). In all of the above cases, we derive from the arguments  in \cite{SeqRPF} various types of limit theorems for non-stationary sequences of random variables of the form $\sum_{j=0}^{n-1}u_j(T_0^nx)$, where $x$ is distributed according to a special Gibbs measure, $T_0^n=T_{n-1}\circ\cdots\circ T_1\circ T_0$ and $u_j$ is a sequence of bounded H\"older continuous functions (in fact, our method allows us to consider a larger class of functions $u_j$).

\section{Preliminaries and main results}
\subsection{Young towers and countable shifts}
We begin with describing the Tower structure from \cite{Y2}.
Let $(\Del_0,\cF_0,m_0)$ be a probability space, $\{\Lam_0^j:\,j\geq1\}$ be a partition of $\Del_0$ (\text{mod} $m_0$), and $R:\Del_0\to\bbN$ be a (return time) function which is constant on each one of the $\Lam_0^j$'s. We  identify each element $x$ in $\Del_0$ with the pair $(x,0)$, and
for each nonnegative integer $\ell$ let the $\ell$-th floor of the tower be defined by
\[
\Del_\ell=\{(x,\ell)\in\Gam_0\times\{\ell\}:\,\,R(x)>\ell\}
\]
and for each $j$ so that $R|\Lam_0^j>\ell$ set 
\[
\Lam_\ell^j=\{(x,\ell)\in\Del_\ell:\,x\in\Lam_0^j\}\subset\Del_\ell.
\]
The tower is defined by 
\[
\Del=\{(x,\ell):\,\ell\geq0,\,\,(x,\ell)\in\Del_\ell\}.
\]
Let $f_0:\Del_0\to\Del_0$ be so that for each $j$ the map $f_0|\Del_0^j:\Del_0^j\to\Del_0$ is bijective (\text{mod} ${m_0}$). The dynamics on the tower is given by the map $F:\Del\to\Del$ defined by
\[
F(x,\ell)=\begin{cases}(x,\ell+1) & \text{ if }R(x)>\ell+1 \\(f_0(x),0)& \text{ if }R(x)=\ell+1\end{cases}.
\]
We think of $(f_0(x),0)$ as the return (to the base $\Del_0$) function corresponding to $F$, and when $R(x)=\ell+1$ we will also write $F^R(x,0):=F(x,\ell)=(f_0(x),0)$. It will also be convenient to set 
$F^R(x,\ell)=F^R(x,0)$ for any $\ell\geq1$ and $(x,\ell)\in\Del_\ell$.
We note that in applications usually $\Del_0$ is a subset of a larger set, and $f_0=f^R$ is the return time function (to $\Del_0$) of a different function $f$ (so that the tower is constructed in order to study statistical properties of $f$). We assume here that the partition $\cR=\{\Lam_\ell^j\}$ is generating in the sense that 
\[
\bigvee_{i=0}^\infty F^{-i}\cR
\]
is a partition into points. For each $k\geq0$ and $x\in\Del$, we will denote the element of the partition 
\[
\bigvee_{i=0}^k F^{-i}\cR
\]
containing $x$ by $\cR_k(x)$ (so that $\{x\}=\cap_{k\geq0}\cR_k(x)$).
\begin{remark}\label{Y1 remark}
The tower structure in \cite{Y1} was defined differently- it was built over possibly invertible maps $f$ and the return time function $R$ satisfied different conditions, but the setup described above fits the projection $\bar F$ of the tower from \cite{Y1} on the  quotient space generated by sliding along stable manifolds (see Section 3.1 in \cite{Y1}). Using the projected tower, it is possible to derive statistical limit theorems for the original tower (and hence for the original map over which the tower is constructed, see Section 4.1 in \cite{Y1}). This is clearly also true in the (sequential) case of partial sums of the form $\sum_{j=0}^{n-1}u_j\circ F^j$ considered in what follows (as we will consider sequences of functions $\{u_j\}$ which are uniformly H\"older continuous).
\end{remark}

Next, we lift the $\sig$-algebra $\cF_0$ to $\Del$ by identifying $\Lam_\ell^j$ with $\Lam_0^j$ and lift the probability measure $m_0$ to a measure on $\Del$, by assigning the mass $m_0(\Gam)$ to each subset $\Gamma$ of each $\Lam_\ell^j$, for any $\ell$ and $j$ so that $R|\Lam_0^j>\ell$. Let us denote the above $\sig$-algebra and measure on $\Del$ also by $\cF_0$ and $m_0$, respectively.
We will always assume that $\int Rdm_0<\infty$ which means that $m_0(\Del)<\infty$. Henceforth we will assume that $m_0$ has been normalized so that $m_0(\Del)=1$.
We will assume here  the tower has exponential tails:
\begin{assumption}\label{Tails ass}
The exist constants $q,p>0$ so that for each $n\geq1$,
\[
m_0\{x:\,R(x)>n\}\leq qe^{-pn}.
\]
\end{assumption}

The (separation) distance on the space $\Del$ is defined as follows: for 
any $x=(x^0,0)$ and $y=(y^0,0)$ in $\Del_0$, denote by $s(x,y)$ the greatest integer $n$ so that $(F^R)^p(x)=f_0^p(x^0)$ and $(F^R)^p(y)=f_0^p(y^0)$ lie in the same $\Lam_0^j$, for all $p\leq n$. If $x=(x^0,\ell)$ and $y=(y^0,\ell)$ belong to the same floor $\Del_\ell$ for some $\ell\geq1$ we  set $s(x,y)=s(x^0,y^0)$. When $x$ and $y$ are not in the same floor we set $s(x,y)=0$. Let $\be\in(0,1)$ and define the distance $d(\cdot,\cdot)$ on $\Del$ by $d(x,y)=\beta^{s(x,y)}$. We will also assume that 
\[
F^R:\Lam_0^j\to\Del_0
\]
and its inverse are both non-singular with respect to $m_0$, and that the Jacobian $JF^R$ is locally Lipschitz continuous in the sense that for any $j\geq1$ and $x,y\in\Del_0^j$,
\begin{equation}\label{Jack Reg}
\left|\frac{JF^R(x)}{JF^R(y)}-1\right|\leq Cd(F^R(x),F^R(y))
\end{equation}
for some constant $C$ which does not depend on $j$. Let the transfer operator $L_0$ be defined by 
\[
L_0f(x)=\sum_{y\in F^{-1}\{x\}}JF(y)^{-1}f(y)
\]
where $J_F$ is the Jacobian of $F$.
Note that on $\Del_\ell,\,\ell>0$ we have $L_0f(x,\ell)=f(x,\ell-1)$, while on $\Del_0$ the members of the set $F^{-1}\{x\}$ are of the form $y=(y^0,\ell)$ with $R(y^0)=\ell+1$, and then $JF(y)=JF^R(y^0,0)$.

Next, for each function $f:\Del\to\bbC$, let  $\|f\|_\infty$ denote its supremum and let $L(f)$ denote the infimum of all possible values $L$ so that for any $\ell$ and $x,y\in\Del_\ell$ we have
\[
|f(x)-f(y)|\leq Ld(x,y).
\]
We will say that $f$ is locally Lipschitz continuous if $\|f\|:=\max\{\|f\|,L(f)\}<\infty$, and let us denote the Banach spaces of all complex valued functions $f$ so that $\|f\|<\infty$ by $\cH$.
We will also assume here that the greatest common divisor of the $R_i$'s equals $1$. In this case, by Theorem 1 in \cite{Y2}, there exists a locally Lipschitz continuous function $h_0$ which is bounded, positive and uniformly bounded away from $0$ so that $L_0h_0=h_0$, $m_0(h_0)=1$ and the measure $\mu=h_0dm_0$ is $F$-invariant and the measure preserving system $(\Del,\mu,\cF_0,F)$ is mixing. Now, for each $\ell\geq0$ set $v_\ell=e^{\frac12\ell p}$ (where $p$ comes from Assumption \ref{Tails ass}). We view $\{v_\ell\}$ as a function $v:\Del\to\bbR$ so that $v|\Del_\ell\equiv v_\ell$, and  we set $m=vm_0$ (which is finite in view of Assumption \ref{Tails ass}) and $h=\frac{h_0}v$. Following \cite{Viv2}, consider the transfer operator $L$ given by 
\[
Lf=\frac{L_0(fv)}{v}.
\]
Then $Lh=h$ and $L^*m=m$ (since $L_0^*m_0=m_0$), and the space $\cH$ is $L$-invariant. In fact (see Lemmas 1.4 and 3.4 in \cite{Viv2}), the operator norms $\|L^n\|$ are uniformly bounded in $n$.

Next, let $u_j:\Del\to\bbR$ be a two sided sequence of locally Lipschitz continuous functions so that $B=\sup_{j\in\bbZ}\|u_j\|<\infty$. For each $j$ and $z\in\bbC$ let the transfer operator $\cL_z^{(j)}$ be defined by 
\[
\cL_z^{(j)}f=L(fe^{zu_j}).
\]
Then, for each integer $j$ and $z\in\bbC$, the space $\cH$ is $\cL_z^{(j)}$- invariant (since $e^{zu_j}$  are members of $\cH$ ).
Since the map $z\to e^{zu}\in\cH$ is analytic, the operators $\cL_z^{(j)}$ are analytic functions of in $z$, when viewed as  maps to the space of continuous linear operators $A:\cH\to\cH$, equipped with the operator norm. For each integer $j$, a complex number $z$ and $n\in\bbN$ set 
\[
S_{j,n}u=\sum_{k=0}^{n-1}u_{j+k}\circ F^k
\]
and 
\[
\cL_z^{j,n}=\cL_z^{(j+n-1)}\circ\cdots\circ\cL_z^{(j+1)}\circ \cL_z^{(j)}
\]
which satisfy $\cL_z^{j,k}g=\cL_0^{j,k}(ge^{zS_{j,n}u})=L^n(ge^{zS_{j,n}u})$. When it is more convenient, we will also write 
$L^n=\cL_0^n=\cL_0^{j,n}$.
Henceforth, we will refer to $q,p$ from Assumption \ref{Tails ass}, $C$ from (\ref{Jack Reg}) and  $B=\sup_{j}\|u_j\|$  as the ''initial parameters". Our main result is the following

\begin{theorem}\label{RPF}
There exists a constant $r>0$, which depends only on the initial parameters, so that for any $z\in U= B(0,r):=\{\zeta\in\bbZ:\,|\zeta|<r\}$
there exist families
$\{\la_j(z):\,j\in\bbZ\}$, $\{h_j^{(z)}:\,j\in\bbZ\}$ and $\{\nu_j^{(z)}:\,j\in\bbZ\}$ 
consisting of a nonzero complex number 
$\la_j(z)$, a complex function $h_j^{(z)}\in\cH$ and a 
complex continuous linear functional $\nu_j^{(z)}\in\cH^*$ such that:

(i) For any $j\in\bbZ$,
 $\la_j(0)=1$, $h_j^{(0)}=h$, $\nu_j^{(0)}=m$
 and for any $z\in B(0,r)$,
\begin{equation}\label{RPF deter equations-General}
\cL_z^{(j)} h_j^{(z)}=\la_j(z)h_{j+1}^{(z)},\,\,
(\cL_z^{(j)})^*\nu_{j+1}^{(z)}=\la_j(z)\nu_{j}^{(z)}\text{ and }\,
\nu_j^{(z)}(h_j^{(z)})=\nu_j^{(z)}(h)=1.
\end{equation} 
When $z=t\in\bbR$ then 
 $\la_j(t)>a$ and the function $h_j(t)$ takes values at some interval $[c,d]$, where $a>0$ and $0<c<d<\infty$ depend only on the initial parameters. Moreover, $\nu_j^{(t)}$ is a probability measure which assigns positive mass to open subsets of $\Del$ and the equality 
$\nu_{j+1}(t)\big(\cL_t^{(j)} g)=\la_j(t)\nu_{j}^{(t)}(g)$ holds true for any 
bounded Borel function $g:\Del\to\bbC$.

(ii) Set $U=B(0,r)$. Then the maps 
\[
\la_j(\cdot):U\to\bbC,\,\, h_j^{(\cdot)}:U\to \cH\,\,\text{ and }
\nu_j^{(\cdot)}:U\to \cH^*
\]
are analytic and there exists a constant $C>0$, which depends only on the initial parameters
such that 
\begin{equation}\label{UnifBound}
\max\big(\sup_{z\in U}|\la_j(z)|,\, 
\sup_{z\in U}\|h_j^{(z)}\|,\, \sup_{z\in U}
\|\nu^{(z)}_j\|\big)\leq C,
\end{equation}
where $\|\nu\|$ is the 
operator norm of a linear functional $\nu:\cH\to\bbC$.
Moreover, there exist a constant $c>0$, which depends only on the initial parameters, so that $|\la_j(z)|\geq c$ 
and $\min_{x\in \Del}|h_j^{(z)}(x)|\geq c$ for any integer $j$ and $z\in U$.

(iii) There exist constants $A>0$ and $\del\in(0,1)$, 
which depend only on the initial parameters,
 so that for any $j\in\bbZ$, $g\in\cH$
and $n\geq1$,
\begin{equation}\label{Exponential convergence}
\Big\|\frac{\cL_z^{j,n}g}{\la_{j,n}(z)}-\nu_j^{(z)}(g)h^{(z)}_{j+n}\Big\|\leq A\|g\|\del^n
\end{equation}
where $\la_{j,n}(z)=\la_{j}(z)\cdot\la_{j+1}(z)\cdots\la_{j+n-1}(z)$. 
\end{theorem}

\begin{remark}
It will be clear from the proof of Theorem \ref{RPF} that it also holds true for the transfer operators generated by the countable subshifts of finite type and the potentials $\Phi$ considered in Theorem 1.3 from \cite{Viv1}, with the exception that $|h_j^{(z)}|$ might not be bounded from below (for this we need that $h=h_{j}^{(0)}$ will be bounded from below, which was not obtained in \cite{Viv1}). 
The proof and the statements in this countable subshift case are almost identical to the Young tower case considered above, and so, in order to avoid repetitiveness, we only present the results in the Young tower case. Usually, when applying the RPF theorem in order to obtain limit theorems, the transfer operator is normalized so that the (invaraiant) Gibbs measure $\mu=hd\nu$ becomes a conformal measure. This normalization is not possible when $h$ (or $h_j^{(0)}$) is not bounded away from $0$,  since $h$ appears in the denominator. Still, in the situation in \cite{Viv1} we assume that a conformal measure exists, and so there is no need to normalize the transfer operator, and the result in this paper imply that all the limit theorems in \cite{SeqRPF} hold true for sequences of random variables having the form $\sum_{k=0}^{n-1}u_k(\sig^kx)$, where $\sig$ is the countable shift from \cite{Viv1} and $\mu$ is the conformal measure (which is not necessarily invariant).  
\end{remark}

\subsection{Covering maps}
Let $(\cE_j,d_j)$ be a sequence of metric spaces, which are assumed to be normalized so that $\text{diam}(\cE_j)\leq1$ for each $j$. Let  $T_j:\cE_j\to\cE_{j+1}$ be a sequence of surjective maps. For any $j\in\bbZ$ and $n\in\bbN$ set 
\[
T_j^n=T_{j+n-1}\circ\cdots\circ T_{j+1}\circ T_j.
\] 
\begin{assumption}\label{CovAss}
There exists a constant $\xi>0$ so that for any $j$ and $x,x'\in\cE_{j+1}$ so that $d_{j+1}(x,x')<\xi$ we can write
\[
T_j^{-1}\{x\}=\{y_k:\,k<D\}\,\,\text{ and }\,\,
T_j^{-1}\{x'\}=\{y'_k:\,k<D\}
\]
where $D=D_j(x)\in (1,\infty]$ and for each  $k<D$,
\[
d_j(y_k,y'_k)\leq d_{j}(x,x').
\]
Moreover, there exist constants $\gam>1$ and $n_0\in\bbN$ so that for any $j$ and 
$x,x'\in\cE_{j+n_0}$ with $d_{j+n_0}(x,x')<\xi$ we can write
\[
(T_j^{n_0})^{-1}\{x\}=\{y_k:\,k<D_{n_0}\}\,\,\text{ and }\,\,
(T_j^{n_0})^{-1}\{x'\}=\{y'_k:\,k<D_{n_0}\}
\]
where $D_{n_0}=D_{n_0,j}(x)\in (1,\infty]$ and for each  $k<D_{n_0}$,
\[
d_j(y_k,y'_k)\leq \gam^{-1}d_{j+n_0}(x,x').
\]
\end{assumption}
Assumption \ref{CovAss} holds true for countable non-stationary subshifts of finite type with $n_0=1$. When $\xi>1$ it 
 holds true for families of uniformly distance expanding maps (with $n_0=1$), and in particular for the sequential dynamical systems consider at Section 3 of \cite{Conze}. It also holds true with $\xi>1$
when each $T_j$ is the Gauss map $T:(0,1)\setminus\bbQ\to(0,1)\setminus\bbQ$ given by $T(x)=[x^{-1}]-x^{-1}$ 
(we can take $n_0=2$), as well as in
the setup of non-stationary full-shifts (see Section \ref{examples} for description of these examples). Note that it is possible to obtain our results in the Gauss map case by passing to its symbolic representation (which is a countable full-shift), in which we can take $n_0=1$. Still, our method works fine without using this representation and also works in more general cases in which Assumption \ref{CovAss} holds true only for $n_0>1$, where a symbolic representation might not exist.

Under Assumption \ref{CovAss}, for any $n\geq1$, integers $j$ and $n\geq1$ and $x,x'\in\cE_{j+n}$ with $d_{j+n}(x,x')<\xi$ we can write 
\begin{equation}\label{Large pair}
(T_j^{n})^{-1}\{x\}=\{y_{k,n}:\,k<D_{n}\}\,\,\text{ and }\,\,
(T_j^{n})^{-1}\{x'\}=\{y'_{k,n}:\,k<D_{n}\}
\end{equation}
where $D_n=D_{j,n}(x)\in[1,\infty]$ and for each $k$ we have 
\begin{equation}\label{Iterated distance}
d_j(y_{k,n},y'_{k,n})\leq \gam^{-[\frac{k}{n_0}]}d_{j+n}(x,x').
\end{equation}
For any constants $\al\in(0,1]$ and $Q>0$  let $H_{Q,\al}$ be the collection of all families $\{g_j\}$ of functions $g_j:\cE_j\to\bbR$ so that
for any $k$, $j$ and $x,x'\in\cE_{j+n}$ with $d_{j+n}(x,x')<\xi$ we have
\begin{equation}\label{Regularity of the functions}
\left|\sum_{m=0}^{n-1}\big(g_{j+m}( T_j^my_{k,n})-g_{j+m}( T_j^my'_{k,n})\big)\right|\leq Q\big(d_{j+n}(x,x')\big)^\al
\end{equation}
where the pairs $(y_{k,n},y'_{k,n})$ are the ones satisfying (\ref{Large pair}) and (\ref{Iterated distance}).
This condition is satisfied when all of the $g_j$'s are locally H\"older continuous, uniformly in $j$, but it is also holds true when each $T_j$ is the Gauss map and $g(y)=-2\ln y$, which allows us to consider the classical dynamical system related to continued fractions (see Section \ref{examples}).
For any integer $j$, a number $\al\in(0,1]$ and  $g:\cE_j\to\bbC$ set $\sup|g|=\|g\|_\infty$, 
\[
v_{\al,\xi}(g)=\sup\left\{\frac{|g(x)-g(y)|}{\big(d(x,y)\big)^\al}:\,\,d_j(x,y)\in(0,\xi)\right\}
\]
and $\|g\|_\al=\|g\|_\infty+v_{\al,\xi}(g)$. We denote here by $\cH_j^\al$ the Banach space of all functions $g:\cE_j\to\bbC$ so that $\|g\|_\al<\infty$.

Next, fix some $\al$ and $Q$ and let $\{f_j\},\{u_j\}\in H_{Q,\al}$ be that $B=\sup_{j}\|u_j\|_\infty<\infty$. For each $x\in\cE_{j+1}$ and $z\in\bbC$ write 
\[
\cL_z^{(j)}g(x)=\sum_{y\in T_j^{-1}\{x\}}e^{f_j(y)+zu_j(y)}g(y).
\]
We assume that the  $e^{f_j}$ is are uniformly summable in $j$ the sense that 
\[
\sup_{x\in\cE_j}\sum_{y\in T_j^{-1}\{x\}}e^{f_j(y)}<C
\]
for some constant $C>0$ which does not depend on $j$.
This, in particular, means that $\cL_z^{(j)}g$ are well defined for bounded functions $g$.
For any integer $j$ and $n\geq1$ set $S_{j,n}u=\sum_{k=0}^{n-1}u_{j+k}\circ T_j^k$ and
\[
\cL_z^{j,n}=\cL_z^{(j+n-1)}\circ\cdots\circ\cL_z^{(j+1)}\circ \cL_z^{(j)}.
\]
The following sequential Lasota-Yorke type inequality is  proved in Section \ref{CovSec}:
\begin{proposition}\label{LY covering}
For any $n\in\bbN$, $j\in\bbZ$, $z\in\bbC$
and $g\in \cH_j$,
\begin{eqnarray*}
v_{\al,\xi}(\cL_z^{j,n}g)\leq \|\cL_0^{j,n}\textbf{1}\|_\infty
e^{|\Re(z)|\|S_{j,n} u\|_\infty}\\
\times\big(v_{\al,\xi}(g)(\gam^{-\al[n/n_0]}
+2Q(1+\|z\|_1)\|g\|_\infty\big)
\end{eqnarray*}
where $\|z\|_1=|\Re(z)|+|\Im(z)|$ and $\Re(z)$ ($\Im(z)$) is the 
real (imaginary) part of $z$.  
As a consequence,
\begin{eqnarray}\label{L.Y.-general}
\|\cL_z^{j,n}g\|_{\al,\xi}\leq \|\cL_0^{j,n}\textbf{1}\|_\infty
e^{|\Re(z)|\|S_{j,n} u\|_\infty}\\
\times\big(v_{\al,\xi}(g)(\gam^{-\al[n/n_0]}
+(1+2Q)(1+\|z\|_1)\|g\|_\infty\big).\nonumber
\end{eqnarray}
In particular, 
$\cL_z^{j,n}|_{B_j}:\cH_j\to \cH_{j+n}$ is a continuous linear map.
\end{proposition}


When
Assumption \ref{CovAss} holds true with $\xi\leq1$, we will also need the following
\begin{assumption}\label{RPF 0 ASS}
(i) There exists an integer $m_0$ so that for any $j\in\bbZ$ and  $y\in\cE_j$ we have 
\begin{equation}\label{TopEx}
T_j^{m_0}(B_j(y,\xi))=\cE_{j+m_0}
\end{equation}
where $B_j(y,\xi)$ is an open ball in $\cE_j$ around $y$ with radius $\xi$.
Moreover, the functions $\{f_j\}$ are uniformly bounded (as  $j$ varies). 

(ii) There exists a sequence of strictly positive functions $h_j:\cE_j\to\bbC$ which are  bounded and bounded away from $0$ uniformly in $j$, 
Borelian probability measures $\nu_j$ on $\cE_j$ and positive numbers $\la_j$ which are uniformly bounded in $j$
 so that 
\[
\lim_{n\to\infty}\sup_{j\in\bbZ}\left\|\frac{\cL^{j,n}_0}{\la_{j,n}}-\nu_j\otimes h_{j+n}\right\|_\al=0
\]
where $\la_{j,n}=\prod_{m=0}^{n-1}\la_{j+m}$ and $\nu_j\otimes h_{j+n}:\cH_j\to\cH_{j+n}$ is given by 
$\big(\nu_j\otimes h_{j+n}\big)(g)=\nu_j(g)\cdot h_{j+n}$.
\end{assumption}
When all the $T_j$'s coincide with  the same  countable topologically mixing subshift of finite type $T$ and $f_j=f$ does not depend on $j$, then 
Assumption \ref{RPF 0 ASS} holds true  when $f$ is positive reccurent  with finite pressure and either $T$ has finitely many images or $T$ satisfies the big image condition and $f$ has an invariant Gibbs measure (see Theorem 5 and Proposition 2 in Section 5 of \cite{Sarig1}, Theorem 8 in Section 8 of \cite{Sarig1} and Theorem 1 in \cite{Sarig2}). This assumption also holds true for certain Gibbs-Markov maps (see the corollary proceeding Proposition 1.2 in \cite{Jon}), assuming that (\ref{TopEx}) holds true.



Note that Assumption \ref{RPF 0 ASS} is somehow consistent with the case of Young towers and countable shifts considered in the previous section (i.e. the ones from \cite{Viv1} and \cite{Viv2})- we rely on some {\em a priori} convergence of the system towards an equilibrium. 
Next, we will call the parameters $Q$ and $B$ the initial parameter. Our main result here shows that Theorem 2.4 from \cite{SeqRPF} holds true in the countable degree case (under some conditions):

\begin{theorem}\label{RPF Covering maps}
Suppose that Assumption \ref{CovAss} holds true with $\xi>1$ or that it holds true with $\xi\leq1$ and that, in addition, Assumption \ref{RPF 0 ASS} holds true. Then all the results stated in Theorem \ref{RPF} hold true with triplets that satisfy 
$\la_j(0)=\la_j$, $h_j^{(0)}=h_j$, $\nu_j^{(0)}=\nu_j$ and $\nu_j^{(z)}(\textbf{1})=1$ (instead  of $\nu_j^{(z)}(h_j^{(0)})=1$).
\end{theorem}

\section{Real and complex cones and the associated Hilbert projective metrics: summary}


\subsection{Real cones}
Let $X$ be a real vector space. A subset $\cC_\bbR\subset X$ is called 
a proper real convex cone (or, in short, a real cone if $\cC$ is convex,
invariant under multiplication of nonnegative  numbers and 
$\cC_\bbR\cap-\cC_\bbR=\{0\}$. Next, assume that $X$ is a Banach space
and let $\cC_\bbR\subset X$ be
a closed real cone. For any nonzero elements $f,g$ of $\cC_\bbR$ set 
\begin{equation}\label{ReHilMet beta.0}
\be_{\cC_\bbR}(f,g)=\inf\{t>0: tf-g\in\cC_\bbR\}
\end{equation}
where we use the convention $\inf\emptyset=\infty$. We note that $\be_{\cC_\bbR}(f,g)>0$
since otherwise $-g$ lays in (the closure of) $\cC_\bbR$ which together with the inclusion
$g\in\cC_\bbR$ implies that $g=0$. 
The real Hilbert (projective) metric
$d_{\cC_\bbR}:\cC_\bbR\times\cC_\bbR\to[0,\infty]$
associated with the cone is given by 
\begin{equation}\label{General real Hilbert metric.0}
d_{\cC_\bbR}(f,g)=\ln \big(\be_{\cC_\bbR}(f,g)\be_{\cC_\bbR}(g,f)\big)
\end{equation}
where we use the convention $\ln\infty=\infty$.

Next, let $X_1$ and $X_2$ be two real Banach spaces and let $\cC_i\subset X_i, i=1,2$
be two closed real cones. Let $A:X_1\to X_2$ be a continuous linear
transformation such that $A\cC_1\setminus\{0\}\subset\cC_2\setminus\{0\}$ and set 
\[
D=\sup_{x_1,x_2\in \cC_1\setminus\{0\}}d_{\cC_2}(Ax_1,Ax_2).
\]
The following theorem is a particular case of Theorem 1.1  in \cite{Liver} (see \cite{Bir} for the case when $\cC_1=\cC_2$).
\begin{theorem}
For any nonzero $x,x'\in\cC_1$ we have
\[
d_{\cC_2}(Ax,Ax')\leq\tanh\big(\frac14D\big)d_{\cC_1}(x,x')
\]
where $\tanh \infty:=1$
\end{theorem}
This lemma means that any linear map between two (punctured) real closed cones 
weakly contracts the corresponding Hilbert metrics, and this contraction is strong 
if the (Hilbert) diameter of the image is finite. 

\subsection{Complex cones}
In this section we will present briefly the theory of complex Hilbert metrics developed in \cite{Rug} and 
\cite{Dub1}. We refer the readers' to Appendix A in \cite{book} for a more detailed description of this theory.

Let $Y$ be a complex Banach space. We recall the following definitions
A subset $\cC\subset Y$ is called a \textbf{complex cone} if 
$\bbC'\cC\subset\cC$, where $\bbC'=\bbC\setminus\{0\}$.
The cone $\cC$ is said to be \emph{proper}
 if its closure $\bar{\cC}$ does not contain any two dimensional 
complex subspaces. The \emph{dual cone\em}
 $\cC^*\subset Y^*$ is the set given by
\[
\cC^*=\{\mu\in Y^*:\, \mu(c)\not=0\,\,\,\,\forall c\in\cC'\}
\]
where
$
\cC'=\cC\setminus\{0\},
$
and $Y^*$ is the space of all continuous linear functionals $\mu:Y\to\bbC$
equipped with the operator norm. It is clear that 
$\cC^*$ is a complex cone.
We will say that $\cC$ is \emph{linearly convex}
if for any $x\not\in\cC$
there exists $\mu\in\cC^*$ such that $\mu(x)=0$, i.e. the complement
of $\cC$ is the union of the kernels $Ker(\mu),\mu\in\cC^*$.
Then the complex cone $\cC^*$ is linearly convex, since for any $c\in\cC$
the corresponding evaluation map $\nu\to\nu(c)$ is a member 
of the dual cone $(\cC^*)^*$ of $\cC^*$.

We  introduce now the notation of the complex Hilbert 
projective metric $\del_{\cC}$
of a proper complex cone defined in \cite{Dub1}. Let $x,y\in\cC'$ and consider the set
 $E_{\cC}(x,y)$ given by  
\[
E_{\cC}(x,y)=\{z\in\bbC:\,zx-y\not\in\cC\}.
\]  
Since $\cC$ is proper (and $\bbC'$ invariant) the set $E_{\cC}(x,y)$ is nonempty. 
When $x$ and $y$ are collinear set $\del_\bbC(x,y)=0$ and otherwise set
\[
\del_{\cC}(x,y)=\ln\big(\frac ba\big)\in[0,\infty]
\]
where 
\[
a=\inf|E_\bbC(x,y)|\in[0,\infty]\,\,\,\text{ and }\,\,\,b=\sup|E_\bbC(x,y)|\in[0,\infty]
\]
are the ``largest" and ``smallest" modulus of the set $E_{\cC}(x,y)$, respectively.
Observe that  $\del_{\cC}(x,y)=\del_{\cC}(c_1x,c_2y)$ 
for any $c_1,c_2\in\bbC'=\bbC\setminus\{0\}$, i.e. $\del_\cC$ is projective.
When $\cC$ is linearly convex then $\del_{\cC}$ satisfies the 
triangle inequality and so it is a projective metric (see \cite{Dub1} and \cite{Dub2}). 
We remark that a different notion of a complex Hilbert metric $d_\cC$ was defined in \cite{Rug}, which was prior
to the definition of $\del_\cC$. For canonical complexifications of real cones (defined below)  
$d_{\cC}$ and $\del_{\cC}$ are equivalent (see Section 5 in \cite{Dub1}), and so
it will makes no difference whether we use $\del_{\cC}$ or $d_{\cC}$.

Next, let $X$ be a real Banach and let $\cC_\bbR\subset X$ be  
a real cone.
Let $Y=X_\bbC=X+iX$ be its complexification (see Section 5 from \cite{Rug}).
Following \cite{Rug}, we define the  canonical complexification of $\cC_\bbR$ 
by 
\begin{equation}\label{Complexification}
\cC_\bbC=\{x\in X_\bbC:\,\Re\big(\overline{\mu(x)}\nu(x)\big)
\geq0\,\,\,\,\forall\mu,\nu\in\cC_\bbR^*\}
\end{equation}
where 
\[
\cC^*_\bbR=\{\mu\in X^*:\,\mu(c)\geq0\,\,\,\,\,\forall c\in\cC_\bbR\}
\]
and $X^*$ is the space of all continuous linear functions $\mu:X\to\bbR$ equipped 
with the operator norm.
Then $\cC_\bbC$ is a proper complex cone
(see Theorem 5.5 of  \cite{Rug}) and by
\cite{Rug} and \cite{Dub1}  we have the following 
polarization identities
\begin{equation}\label{polar}
\cC_\bbC=\bbC'\big(\cC_\bbR+i\cC_\bbR\big)=\bbC'
\{x+iy: x\pm y\in\cC_\bbR\}
\end{equation}
where we recall that $\bbC'=\bbC\setminus\{0\}$. 
Moreover, when 
\[
\cC_\bbR=\{x\in X: \mu(x)\geq0\,\,\,\,\,\,\forall \mu\in \cS\}
\]
for some $\cS\subset X^*$, then
 \begin{equation}\label{complexification}
\cC_\bbC=\{x\in X_\bbC:\,\Re\big(\overline{\mu(x)}\nu(x)\big)\geq0\,\,\,\,\,\forall\mu,
\nu\in \cS\}
\end{equation}
since $\cS$ generates the dual cone $\cC_\bbR^*$. Note that by Lemma 4.1 in \cite{Dub2}, 
a canonical complexification $\cC_\bbC$ 
of a real cone $\cC_\bbR$ is linearly convex if there exists
a continuous linear functional which is strictly positive on $\cC_\bbR'=\cC_\bbR\setminus\{0\}$. 

Next, recall that a real cone  $\cC_\bbR$ is said to have a bounded aperture, if there exists $K>0$ and  
a  continuous linear functional $\mu\in\cC_\bbR^*$ so that
\begin{equation}\label{aper1}
\|x\|\|\mu\|\leq K\mu(x)\,\,\text{ for any }\,\, x\in\cC_\bbR.
\end{equation}
Now, following \cite{Rug},  a complex cone 
$\cK_\bbC$ is said to have a bounded aperture if there exists $K>0$ and 
 a continuous  linear functional $\mu\in\cK_\bbC^*$
such that for any $x\in\cK_\bbC$,
\begin{equation}\label{aper2}
\|x\|\|\mu\|\leq K|\mu(x)|.
\end{equation}
The following result appears in \cite{Rug} as Lemma 5.3:
\begin{lemma}\label{Lemma 5.3}
Let $X$ be a real Banach space and let $Y$ be its complexification. 
Let $\cC_\bbR\subset X$ be a real cone, and assume that  
(\ref{aper1}) holds true with some $\mu$ and $K$. Then the complexification $\cC_\bbC$
of $\cC_\bbR$ satisfies (\ref{aper2}) with $\mu_\bbC$ and $2\sqrt 2K$, where
$\mu_\bbC$ is the unique extension of $\mu$ to the complexified space $Y$.
\end{lemma}
We refer the readers to Lemma A.2.7 from Appendix A in \cite{book} for conditions which guarantee that the dual cone of a canonical complexification of a real cone has bounded aperture (see also the beginning of Section \ref{Sec 4.1}).

Next, the following  assertion is formulated as
Theorem 3.1 in \cite{Dub2} and it summarizes some of the 
main results from \cite{Dub1}.  
\begin{theorem}\label{thm3.1}
Let $\big(X,\|\cdot\|\big)$ be a complex Banach spaces and $\cC\subset X$ be a complex cone.

(i)
Suppose that the cone $\cC$  is linearly convex and of
bounded (sectional) aperture. Then $(\cC'/\sim,\del_{\cC})$ is a complete metric
space, where $x\sim y$ if and only if $\bbC'x=\bbC'y$.

(ii) Let $K>0$ 
and $\mu\in\cC^*$ be such that (\ref{aper2}) 
holds true for any $x\in\cC$. Then for any $x,y\in\cC'$,
\[
\Big\|\frac x{\mu(x)}-\frac y{\mu(y)}\Big\|\leq\frac K{2\|\mu\|}\del_{\cC}(x,y).
\]

(iii) Let $\cC_1$  be a complex cone in some complex 
Banach space $X_1$, and  $A:X\to X_1$ be a complex linear map such that 
$A\cC'\subset\cC_1'$. 
Set 
$
\Del=\sup_{u,v\in\cC'}\del_{\cC_1}(Au,Av)
$
and assume that $\Del<\infty$. 
Then for any $x,y\in\cC'$, 
\[
\del_{\cC_1}(Ax,Ay)\leq\tanh\big(\frac\Del4\big)\del_{\cC}(x,y).
\] 
\end{theorem}
Theorem \ref{thm3.1} $(ii)$ means that any linear map between two (punctured) cones
whose image has finite $\del_{\cC_1}$ (Hilbert) diameter
is a weak contraction with respect to the appropriate projective metrics and that this 
contraction is strong.

\section{Cones for Young towers and countable shifts}\label{YT CntblShift}

First, for any $\ve_0>0$ and $s\geq1$ we can partition $\Del$ into a finite number of disjoint sets $P_2$ and $P',\,P'\in\cP_1$ so that $m(P_2)<\ve_0$ and the diameter each one of the $P'$'s is less than $\gam_s$, where $\gam_s\to 0$ when $s\to\infty$. One way to construct such partitions is as in \cite{Viv2}, and another way is to take a finite collection $\Gam_s$ of the $\Del_\ell^j$'s so that the set 
\[
P_2=\bigcup_{i=0}^s\big(F^R\big)^{-i}\bigcup_{\Del_\ell^j\not\in\Gam_s}\Del_\ell^j
\]
satisfies $m(P_2)<\ve_0$. Denote the above partition by $\cP$. Note that since $\cP$ is finite, then by applying Theorem 1.2 in \cite{Viv2} we deduce that for any $0<\al<1<\al'$ there exists $q_0$ so that for any $k\geq q_0$ and $P,P'\in\cP$ we have 
\[
\al<\frac{m(P\cap F^{-k}P')}{m(P)\mu(P')}<\al'.
\]

Following \cite{Viv2} (and \cite{Viv1}), for any $a,b,c>0$ let the real cone $\cC_{a,b,c,\ve_0,s}$ consists of all the real valued  locally Lipschitz continuous functions $f$ so that:
\begin{itemize}
\item
$0\leq\frac{1}{\mu(P)}\int_P fdm=\frac{1}{\mu(P)}\int_P f/h d\mu\leq a\int fdm;\,\,\forall\,P\in\cP$.
\\
\item
$L(f)\leq b\int fdm$.
\\
\item
$|f(x)|\leq c\int fdm,\,\,\text{for any }\,x\in P_2$
\end{itemize}
If $f\in\cC_\bbR$ then for any $x\in\Del\setminus P_2$,
\[
|f(x)|\leq\frac{1}{m(P_1(x))}\int_{P_1(x)} fdm+\gam_sL(f)\leq (a\|h\|_\infty+b\gam_s)\int fdm
\]
where $P_1(x)\in\cP_1$ is the partition element containing $x$, and we used that $\mu=hdm$. Therefore, with 
\[
c_1=c_1(s,a,b)=a\|h\|_\infty+b\gam_s
\]
and $c_2=\max\{c,c_1\}$ we have 
\begin{equation}\label{f bound}
\|f\|_\infty\leq c_2\int fdm.
\end{equation}
This essentially means that we could have just required that the third condition holds true for any $x\in\Del$, and not only in $P_2$ (by taking $c>c_1$).
Note that if $\int\cL_0^kfdm=0$ for some $k$ and
$f\in\cC_{a,b,c,\ve_0,s}$ then, since 
\[
\int\cL_0^kfdm=\int fdm=0
\]
it follows from (\ref{f bound}) that $f=0$. This means that if, for some $k$, the cone $\cC_{a,b,c,\ve_0,s}$ is $\cL_0^k$-invariant then $\cL_0^k$ is strictly positive with respect to this cone.

The following result was (essentially) proved as in \cite{Viv2}:
\begin{theorem}\label{Viv2Thm}
There exists $\ve_0,s,a,b,c>0$, $\sig\in(0,1)$, $k_0\in\bbN$ and $d_0>0$ so that  with $\cC_\bbR=\cC_{a,b,c,\ve_0,s}$, for any $k\geq k_0$ we have 
\[
L^k\cC\subset\cC_{\sig a,\sig b,\sig c,\ve_0,s}
\]
and for any $f,g\in\cC_\bbR$,
\[
d_{\cC_\bbR}(L^kf,L^kg)\leq d_0.
\]
\end{theorem}
Let us denote by $\cC$ the canonical complexification of the real cone $\cC_\bbR$ from Theorem \ref{Viv2Thm}.
The main result in this section is the following

\begin{theorem}\label{YT cones thm}
For any sufficiently large $a,b,c$ and $d$ we have:

(i) The cone $\cC$ is linearly convex, it contains the functions $h$ and $\textbf{1}$ (the function which takes the constant value $1$). Moreover,  the measure $m$, when viewed as a linear functional, is a member of $\cC_\bbR^*$ and
the cones  $\cC$ and $\cC^*$ have bounded aperture. In fact,
there exist constants $K,M>0$ so that for any $f\in\cC$ and $\mu\in\cC^*$, 
\begin{equation}\label{aperture0}
\|f\|\leq K|m(f)|
\end{equation}
and
\begin{equation}\label{aperture}
\|\mu\|\leq M|\mu(h)|.
\end{equation}

(ii) The cone $\cC$ is reproducing. In fact, there exists a constant $K_1$ so that for any $f\in\cH$ there exists $R(f)\in\bbC$ so that $|R(f)|\leq K_1\|f\|$
and 
\[
f+R(f)h\in\cC.
\]

(iii) There exist  constants $r>0$ and $d_1>0$ so that for any integer $j$, a complex number $z\in B(0,r)$ and $k_0\leq k\leq 2k_0$, where $k_0$ comes from Theorem \ref{Viv2Thm}, we have 
\[
\cL_z^{j,k}\bbC'\subset\bbC'
\]
and 
\[
\sup_{f,g\in\cC'}\del_{\cC}(\cL_z^{j,k}f,\cL_z^{j,k}g)\leq d_1
\]
where $\cC'=\cC\setminus\{0\}$.
\end{theorem}
Theoerm \ref{RPF} follows from Theorem \ref{YT cones thm}, exactly as in Section 3 of \cite{SeqRPF}, and our main task in this Section is to prove Theorem 
\ref{YT cones thm}.

\subsection{Proof of Theorem \ref{YT cones thm}}\label{Sec 4.1}
We begin with the proof of the first part. First, since 
\[
\int_A hdm=\int Ad\mu=\mu(A)
\]
for any measurable  set $A$, 
it is clear that $h\in\bbC_\bbR$ if $a>1$, $b>L(h)$ and $c>\|h\|_\infty$. Moreover, if $c>1$ and $a>D$, where
\begin{equation}\label{D def}
D=\max\Big\{\frac{m(P)}{\mu(P)}:\,P\in\cP\Big\}
\end{equation}
then $\textbf{1}\in\bbC_\bbR$.

Next, if $f\in\bbC_\bbR'$ and $m(f)=0$ then by (\ref{f bound}) we have $f=0$ and so $m\in\cC_\bbR^*$ (since $m\geq0$ on $\cC_\bbR$). In fact, it follows
from the definitions of the norm $\|f\|$ and from (\ref{f bound}) that 
\[
\|f\|\leq\|f\|_\infty+L(f)\leq (c_2+b)m(f)=(c_2+b)\int fdm
\]
and therefore by Lemma \ref{Lemma 5.3} the inequality (\ref{aperture0}) hold true with $K=2\sqrt 2(c_2+b)$. According to Lemma A.2.7 in Appendix A of \cite{book}, for any $M>0$,
inequality (\ref{aperture}) holds true for any $\mu\in\cC^*$ if
\begin{equation}\label{Const M}
\{x\in X: \|x-h\|<\frac1M\}\subset\cC.
\end{equation}
Now we will show how to find a constant $M$ for which (\ref{Const M}) holds true.
For any $f\in\cH$, $P\in\cP$ and $x_1\in P_2$, and distinct $x,y$ which belong to the same level  $\Del_\ell$ (for some $\ell$) set
\begin{eqnarray*}
\Upsilon_P(f)=\frac1{\mu(P)}\int_P fdm,\,\,
\Gam_P(f)=a\int fdm-\frac1{\mu(P)}\int_P fdm,\\
\Gam_{x,y}(f)=b\int fdm-\frac{f(x)-f(y)}{d(x,y)}
\, \text{ and }\,
\Gam_{x_1,\pm}(f)=c\int fdm\pm f(x_1).
\end{eqnarray*}
Let $\cS$ be the collection of all the above linear functionals. Then 
\[
\cC_\bbR=\{f\in\cH:\,s(f)\geq0,\,\forall s\in\cS\}
\]
and so
\begin{equation}\label{Complexification1}
\cC_\bbC=\{f\in \cH:\,\Re\big(\overline{\mu(f)}\nu(f)\big)
\geq0\,\,\,\,\forall\mu,\nu\in\cS\}.
\end{equation}
Let $g\in\cH$ be of the form $g=h+q$ for some $q\in\cH$. We need to find a constant $M>0$ so that $h+q\in \cC$ if $\|q\|<\frac1M$. In view of (\ref{Complexification1}), there are several cases to consider. First, suppose that $\nu=\Upsilon_{P}$ and $\mu=\Upsilon_Q$ for some $P,Q\in\cP$. Since
\[
\frac{1}{\mu(A)}\int hdm=\frac{1}{\mu(A)}\int 1d\mu=1
\]
for any measurable set $A$ with positive measure, we have 
\[
\Re\big(\overline{\mu(h+q)}\nu(h+q)\big)\geq 1-(D^2\|q\|^2+2D\|q\|)\geq 
1-(D+\|q\|)^2
\]
where $D$ was defined in \ref{D def}.
Hence
\[
\Re\big(\overline{\mu(h+q)}\nu(h+q)\big)>0
\]
if $\|q\|$ is sufficiently small. Now consider the case when $\mu=\Upsilon_P$ for some $P\in\cP$ and
$\nu$ is one of the $\Gamma$'s, say $\nu=\Gam_{x,y}$. Then
\begin{eqnarray*}
\Re\big(\overline{\mu(h+q)}\nu(h+q)\big)\geq b-\|h\|-b\|q\|-\|q\|-D\|q\|(b+\|h\|+b\|q\|+\|q\|)\\
\geq b-\|h\|-C(D,b)(\|h\|+\|q\|+\|q\|^2) 
\end{eqnarray*}
where $C(D,b)>0$ depends only on $D$ and $b$. If $\|q\|$ is sufficiently small and $b>\|h\|$ then 
the above left hand side is clearly positive. Similarly, if $\|h\|<\min\{a,b,c\}$ and $\|q\|$ is sufficiently small then 
\[
\Re\big(\overline{\mu(h+q)}\nu(h+q)\big)>0
\]
when either $\nu=\Gam_{x_1,\pm}$ or $\nu=\Gam_{x,y}$.

Next, consider the case when $\mu=\Gam_{x_1,\pm}$ for some $x_1\in P_2$ and  $\nu=\Gam_{x,y}$ for some distinct $x$ and $y$ in the same floor. Then
\begin{eqnarray*}
\Re\big(\overline{\mu(h+q)}\nu(h+q)\big)\geq 
cb-cb\|q\|^2-cb\|q\|\\-c\|q\|(\|h\|+\|q\|)-cb\|q\|-c\|h\|-c\|q\|-(\|q\|+\|h\|)(\|q\|+b+\|h\|+\|q\|)
\end{eqnarray*}
where we used that $\int hdm=1$. Therefore, if $\|q\|$ is sufficiently small and $c$ and $b$ are sufficiently large then
\[
\Re\big(\overline{\mu(h+q)}\nu(h+q)\big)>0.
\] 
Similarly, since 
\[
\left|\frac1{\mu (P)}\int_P qdm\right|\leq D\|q\|
\]
for any other choice of $\mu,\nu\in\cS\setminus\{\Upsilon_P\}$,
we have
\begin{eqnarray*}
\Re\big(\overline{\mu(h+q)}\nu(h+q)\big)\geq \min\{a^2,ab,ac,b^2,bc,c^2\}(1-9\|q\|-9\|q\|^2)\\-9\max\{a,b,c,aD,bD,cD,1\}(\|h\|+\|q\|)^2
\end{eqnarray*}
and so, when $a,b,c$ are sufficiently large and $\|q\|$ is sufficiently small then the above left hand side is positive. The proof of Theorem \ref{YT cones thm} (i) is now complete. 

The proof of Theorem \ref{YT cones thm} (ii) proceeds exactly as the proof of Lemma 3.11 in \cite{Viv2}: for a real valued function $f\in\cH$, it is clearly enough to take any $R(f)>0$ so that
\begin{eqnarray*}
R(f)>(a-1)^{-1}\cdot\max\Big\{\frac1{\mu(P)}\int_P fdm-a\int fdm:\,\,P\in\cP\Big\},\\
R(f)>\frac{L(f)-b\int fdm}{b-L(h)},\,\,R(f)>\max\Big\{-\frac1{\mu(P)}\int_P fdm:\,\,P\in\cP\Big\}\,\,\text{ and }\\
R(f)>\frac{c\int fdm-\|f\|_\infty}{c-\|h\|_\infty}
\end{eqnarray*}
where we take $a,b$ and $c$ so that all the denominators appearing in the above inequalities  are positive,
and we used that $\frac{1}{\mu(A)}\int hdm=1$ for any measurable set $A$ (apply this with $A=P\in\cP$).
For complex valued $f$'s we can write $f=f_1+if_2$, then take $R(f)=R(f_1)+iR(f_2)$ and use that with $\bbC'=\bbC\setminus\{0\}$, 
\[
\cC=\bbC'(\cC_\bbR+i\cC_\bbR).
\]

Now we will prove Theorem \ref{YT cones thm} (iii). Let $k_0\leq k\leq 2k_0$, where $k_0$ comes from Theorem \ref{Viv2Thm}.
According to Theorem A.2.4  in Appendix A of \cite{book} (which is Theorem 4.5 in \cite{Dub2}), if 
\begin{equation}\label{Comp1}
|s(\cL_z^{j,k} f)-s(\cL_0^{k}f)|\leq \ve s(\cL_0^{j,k}f)
\end{equation}
for any nonzero $f\in\cC_\bbR$, for some $\ve>0$ so that 
\[
\del:=2\ve\Big(1+\cosh\big(\frac12 d_0\big)\Big)<1
\]
where $d_0$ comes from Theorem \ref{Viv2Thm},
then, with $\cC'=\cC\setminus\{0\}$,
\begin{equation}\label{I}
\cL_z^{j,k}\cC'\subset\cC'
\end{equation}
and
\begin{equation}\label{II} 
\sup_{f,g\in\cC}(\cL_z^{j,k}f,\cL_z^{j,k}g)\leq d_0+6|\ln(1-\del)|.
\end{equation}
We will show now that there exists a constant $r>0$ so that (\ref{Comp1}) holds true for any $z\in B(0,r)$ and $f\in\cC_\bbR$. We need first the following very elementary result, which for the sake of convenience is formulated here as a lemma.
\begin{lemma}\label{A A' lemma}
Let $A$ and $A'$ be complex numbers, $B$ and $B'$ be real numbers, and let $\ve_1>0$ and $\sig\in(0,1)$ so that
\begin{itemize}
\item
$B>B'$
\item
$|A-B|\leq\ve_1B$
\item
$|A'-B'|\leq\ve_1 B $
\item
$|B'/B|\leq\sigma$.
\end{itemize}
Then 
\[
\left|\frac{A-A'}{B-B'}-1\right|\leq 2\ve_1(1-\sig)^{-1}.
\]
\end{lemma}
The proof of Lemma \ref{A A' lemma} is very elementary, just write
\[
\left|\frac{A-A'}{B-B'}-1\right|\leq\left|\frac{A-B}{B-B'}\right|+
\left|\frac{A'-B'}{B-B'}\right|\leq \frac{2B\ve_1}{B-B'}=\frac{2\ve_1}{1-B'/B}.
\]
Next, let $f\in\cC_\bbR'$. First, suppose that $s$ have the form 
$s=\Gam_{P}$ for some $P\in\cP$. Set
\begin{eqnarray*}
A=a\int \cL_z^{j,k}fdm,\,\, A'=\frac1{\mu(P)}\int_P \cL_z^{j,k}fdm,\\
B=a\int\cL_0^{j,k}fdm\,\,\text{ and }\,\,B'=\frac1{\mu(P)}\int_P \cL_0^{j,k}fdm.
\end{eqnarray*}
Then $B=a\int fdm$ (since $m$ is conformal) and
\[
|s(\cL_z^{j,k})-s(\cL_0^{j,k})|=|A-A'-(B-B')|.
\]
We want to show that the conditions of Lemma \ref{A A' lemma} hold true. 
By Theorem \ref{Viv2Thm} we have 
\begin{equation}\label{Smaller cone}
\cL_0^{j,k}f\in\cC_{\sig a,\sig b,\sig c,s,\ve_0}
\end{equation}
which in particular implies that 
\[
0\leq B'\leq \sig a\int\cL_0^{j,k}fdm=\sig B.
\]
Since $f$ is nonzero and $\int\cL_0^{j,k}fdm=\int fdm\geq 0$ the number $B$ is positive 
(since (\ref{aperture0}) holds true). It follows that $B>B'$ and that 
\[
|B'/B|\leq\sig<1. 
\]
Now we will estimate $|A-B|$. For any complex $z$ so that $|z|\leq1$ write
\begin{eqnarray*}
|A-B|=a\left|\int\cL_0^{k}\big(f(e^{zS_{j,k}u}-1)\big)dm\right|\leq a\|f\|_\infty\|e^{zS_{j,k}u}-1\|_\infty
\int\cL_0^k\textbf{1}dm\\= a\|f\|_\infty\|e^{zS_{j,k}u}-1\|_\infty\int \textbf{1}dm=
a\|f\|_\infty\|e^{zS_{j,k}u}-1\|_\infty\\
\leq ac_2\int fdm\,\cdot(2k_0R\cdot|z|\|u\|_\infty)=2ac_2k_0R\|u\|_\infty|z|\int\cL_0^{k}f dm=R_1|z|B
\end{eqnarray*}
where $\textbf{1}$ is the function which takes the constant value $1$, $\|u\|_\infty=\sup_{j\in\bbZ}\|u_j\|_\infty$, $R$ is some constant which depends only on $k_0$ and $\|u\|_\infty$ 
and 
\[
R_1=2c_2k_0R\|u\|_\infty.
\]
In the latter estimates we have also used (\ref{f bound}).
It follows  that the conditions of Lemma \ref{A A' lemma} are satisfied with $\ve=R_1|z|$. Now we will estimate $|A'-B'|$. 
First, write 
\begin{eqnarray*}
|A'-B'|\leq\frac1{\mu(P)}\int_P\big|\cL_z^{j,k}f-\cL_0^{j,k}f\big|dm=
\frac1{\mu(P)}\int_P\big|\cL_0^{k}\big(f(e^{zS_{j,k}u}-1)\big)|dm\\
\leq\|f\|_\infty \|e^{zS_{j,k}u}-1\|_\infty\frac1{\mu(P)}\int_P\cL_0^k\textbf{1}dm
=\|f\|_\infty \|e^{zS_{j,k}u}-1\|_\infty\frac{m(P)}{\mu(P)}\\\leq
 Dc_2\int fdm\,\cdot 2k_0R\|u\|_\infty|z|
=R_2|z|B
\end{eqnarray*}
where $D$ is defined by (\ref{D def}) and
\[
R_2= Da^{-1}2c_2k_0R\|u\|_\infty.
\]
We conclude now from Lemma \ref{A A' lemma} that 
\[
|s(\cL_z^{j,k})-s(\cL_0^{j,k})|\leq 2R_3(1-\sig)^{-1}|z|s(\cL_0^{j,k})
\]
where $R_3=\max(R_1,R_2)$.

Next, consider the case when $s$ have the form $s=\Gam_{x,\pm}$ for some $x\in\Del$. Set
\begin{eqnarray*}
A=c\int \cL_z^{j,k}fdm,\,\, A'=\pm \cL_z^{j,k}f(x),\\
B=c\int\cL_0^{j,k}fdm\,\,\text{ and }\,\,B'=\pm\cL_0^{j,k}f(x).
\end{eqnarray*}
Then $B>0$ and by (\ref{Smaller cone}) we have 
\[
|B'|\leq \sig B.
\]
Similarly to the previous case, we have 
\[
|A-B|\leq R_4B|z|
\]
where $R_4=2c_2k_0R\|u\|_\infty=2k_0R\|u\|_\infty$. Now we will estimate $|A'-B'|$. Using (\ref{f bound}) 
we have
\begin{eqnarray*}
|A'-B'|=|\cL_z^{j,k}f(x)-\cL_0^{j,k}f(x)|\leq \|f\|_\infty\|e^{zS_{j,k}u}-1\|_\infty\cL_0^{j,k}\textbf{1}(x)\\
\\\leq c_2\int fdm\,\cdot(2k_0|z|R\|u\|_\infty M_1)=BR_5|z| 
\end{eqnarray*}
where $R_5=2c_2k_0|z|R\|u\|_\infty M_1$ and
$M_1$ is an upper bound on the values of $\|L^k\textbf{1}\|_\infty$ for $k_0\leq k\leq 2k_0$ (in fact, we can use Lemma 1.4 in \cite{Viv2} and obtain an upper bound which does not depend on $k_0$).  
Since 
\[
|s(\cL_z^{j,k})-s(\cL_0^{j,k})|=|A-A'-(B-B')|,
\]
we conclude from Lemma \ref{A A' lemma} that 
\[
|s(\cL_z^{j,k})-s(\cL_0^{j,k})|\leq 2R_6(1-\sig)^{-1}|z|s(\cL_0^{j,k})
\]
where $R_6=\max\{R4,R_5\}$. 

Finally, we consider the case when $s=\Gam_{x,x'}$ for some distinct $x'$ and $x'$ which belong to the same floor of $\Del$. 
Set 
\begin{eqnarray*}
A=b\int \cL_z^{j,k}fdm,\,\, A'=\frac{\cL_z^{j,k}f(x)-\cL_z^{j,k}f(x')}{d(x,x')},\\
B=b\int\cL_0^{j,k}fdm\,\,\text{ and }\,\,B'=\frac{\cL_0^{j,k}f(x)-\cL_0^{j,k}f(x')}{d(x,x')}.
\end{eqnarray*}
Then, exactly as in the previous cases, $B>0$,  $|B'|\leq \sig B$,
\[
|s(\cL_z^{j,k})-s(\cL_0^{j,k})|=|A-A'-(B-B')|
\]
and 
\[
|A-B|\leq R_7B|z|
\]
where $R_7=2c_2k_0R\|u\|_\infty$.  Now we will estimate $|A'-B'|$.
Let $\ell$ be so that $x,x'\in\Del_\ell$ and write
$x=(x_0,\ell)$ and $x'=(x_0',\ell)$. Then $d(x,x')=d((x_0,m),(x_0',m))$ for any $0\leq m\leq\ell$.
If $k\leq \ell$ then for any $z$,
\[
\cL_z^{j,k}f(x)=v_\ell^{-1}v_{\ell-k}e^{sS_{j,k}u(x_0,\ell-k)}f(x_0,\ell-k)
\]
and a similar equality hold true with $x'$ in place of $x$.
Set 
\begin{eqnarray*}
U(z)=f(x_0,\ell-k)e^{zS_{j,k}u(x_0,\ell-k)} \,\text{ and }\,  V(z)=f(x_0',\ell-k)e^{zS_{j,k}u(x_0',\ell-k)}
\end{eqnarray*}
and $W(z)=U(z)-V(z)$.
Then for any $z\in\bbC$ so that $|z|\leq1$ we have
\[
d(x,x')|A'-B'|=v_\ell^{-1}v_{\ell-k}|W(z)-W(0)|\leq |z|\sup_{|\zeta|\leq1}|W'(\zeta)|.
\]
Since the functions $u_m,\,m\in\bbZ$ and $f$ are locally Lipschitz continuous (uniformly in $j$) we obtain that for any $\zeta$ so that $|\zeta|\leq1$,
\[
|W'(\zeta)|\leq C_1d(x,x')\|f\| \leq d(x,x')C_1(b+c_2)\int fdm=d(x,x')C_1b^{-1}(b+c_2)B
\]
where $C_1$ depends only on $k_0$ and $\|u\|=\sup_{j\in\bbZ}\|u_j\|$.

Next, suppose that $k>\ell$, where $\ell$ is such that $x,x'\in\Del_\ell$.
The approximation of  $|A'-B'|$ in this case relies on classical arguments from the theory of distance expanding map.
Since $k>\ell$  we can write 
\[
F^{-k}\{x\}=\{y\},\,\,F^{-k}\{x'\}=\{y'\}
\] 
where both sets are at most countable, the map $y\to y'$ is bijective and satisfies that for any $0\leq q\leq k$
\[
d(F^qy,F^qy')\leq \beta^{m_q(y)} d(x,x')\leq d(x,x').
\]
Here $m_q(y)$ is the number of the point among $F^{q+m}y$, $0\leq m\leq k-q$ which belong
to the base $\Del_0$ (so $m_0(y)\geq1$, since $\ell<k$). Note also that the pairs $(y,y')$ also belong to the same partition element $\Del_\ell^j$.
Then for any complex $z$ we have
\[
\cL_z^{j,k}f(x)=v_\ell^{-1}\sum_{y}v(y)JF^k(y)^{-1}e^{zS_{j,k}u(y)}f(y)
\] 
and 
\[
\cL_z^{j,k}f(x')=v_\ell^{-1}\sum_{y}v(y)JF^k(y')^{-1}e^{zS_{j,k}u(y')}f(y')
\] 
where we note that $v(y)=v(y')$ since $y$ and $y'$ belong to the same floor. For any $y$ set 
\[
U_{y}(z)=JF^k(y)^{-1}e^{zS_{j,k}u(y)}f(y)
\]
and 
\[
W_{y,y'}(z)=U_{y}(z)-U_{y'}(z).
\]
Then for any complex $z$ so that $|z|\leq1$ we have
\[
|W_{y,y'}(z)-W_{y,y'}(0)|\leq|z|\sup_{|\zeta|\leq 1}|W'_{y,y'}(\zeta)|.
\]
Since $J_F^{R}$ satisfies (\ref{Jack Reg}) and $u_j$ and $f$ are locally Lipschitz continuous (uniformly in $j$) we 
derive that 
\begin{equation}\label{Der Bound}
\sup_{|\zeta|\leq 1}|W'_{y,y'}(\zeta)|\leq C_2\|f\|d(x,x')(JF^k(y)^{-1}+JF^k(y')^{-1})
\end{equation}
for some constant $C_2$ which depends only on $\|u\|,k_0$ and on $C$ from (\ref{Jack Reg}).
Using that 
\[
\|f\|\leq(c_2+b)\int fdm
\]
we derive now from (\ref{Der Bound})  that 
\begin{eqnarray*}
d(x,x')|A'-B'|=v_\ell^{-1}\left|\sum_{y}v(y)\big(W_{y,y'}(z)-W_{y,y'}(0)\big)\right|
\\\leq \big(|z|d(x,x')C_2\|f\|\big)v_\ell^{-1}\sum_{y}v(y)(JF^k(y)^{-1}+JF^k(y')^{-1})\\=
\big(|z|d(x,x')C_2\|f\|\big)\cdot\big(\cL_0^k\textbf{1}(x)+\cL_0^k\textbf{1}(x')\big)\leq E_1|z| B
\end{eqnarray*}
where $E_1=2M_1C_2b^{-1}(c_2+b)$ and $M_1$ is an upper bound of $\sup_n\|\cL_0^n\|_\infty$.
We conclude that there exists a constant $C_0$ so that for any $s\in\cS$, $f\in\cC'$, $z\in\bbC$ and $k_0\leq k\leq 2k_0$,
\[
|s(\cL_z^{j,k})-s(\cL_0^{j,k})|\leq C_0|z|s(\cL_0^{j,k}).
\]
Let $r>0$ be any positive number so that 
\[
\del_r:=2C_0r\Big(1+\cosh\big(\frac12 d_0\big)\Big)<1.
\]
Then, by (\ref{Comp1}) and what proceeds it, (\ref{I}) and (\ref{II}) hold true
for any $z\in\bbC$ with $|z|<r$, an integer $j$ and $k_0\leq k\leq 2k_0$, and the proof of Theorem \ref{YT cones thm} is complete.
\qed

\begin{remark}
The results from \cite{Viv2} hold true also in the case of summable tails, which do not necessarily decay exponentially fast to $0$. These results were obtained by showing that there exist a sequence of cones $\cC_{j,\bbR},\,j\geq 0$, so that $\cL_0^{k}\cC_{j}\subset\cC_{j+1}$ and the projective diameter $d_j$ of the image of $\cL_0^k\cC_{j,\bbR}$ inside $\cC_{j+1,\bbR}$ converges to infinity with  some rate, which made it possible to apply successively the contraction properties of real Hilbert metric, and obtain a certain type of (not exponentially fast) convergence. When using complex projective metrics (with the canonical complexifications $\cC_j$ of the $\cC_{j,\bbR}$'s),  the magnitude of the perturbation that is allowed, depends on the diameter of the image of the real operators (see (\ref{Comp1})-(\ref{II})). In our situation, this means that it is possible to obtain that $\cL_z^{m,k}\cC_j\subset\cC_{j+1}$ only if 
the term $|z|\|S_{m,k}u\|_\infty$ would be smaller than some $\del_j$, where $\del_j\to\infty$ when $j\to\infty$.
Therefore, it is not possible to choose a neighborhood of $0$ which is independent of $j$ when the tails do not decay exponentially fast.
\end{remark}

\section{Cones for countable covering maps}\label{CovSec}
In this section we will show that the standard cones of locally logarithmic H\"older continuous functions that are commonly used in the setup of finite degree covering maps  (see \cite{MSU}, \cite{Rug}, \cite{Dub2}, \cite{book}, \cite{Kifer-1996}, \cite{Kifer-Thermo} and references therein) can be used successfully in our countable degree setup, as well.

In what follows, when it is more convenient, we will denote the norm  $\|\cdot\|_\al$ also by $\|\cdot\|$ (i.e. we will omit the subscript $\al$). We begin with proving Proposition \ref{LY covering}. The proof proceeds essentially as the proof of Lemma 5.6.1 in \cite{book}, and we include here all the details for readers' convenience.
Let $n,z$ and $g$ be as in the statement of the the proposition.
We begin with approximating $v_{\al,\xi}(\cL_z^{\om,n}g)$.  Let
$x,x'\in\cE_{j+n}$ be such that
$d_{n+j}(x,x')<\xi$ and let $\{y_{i,n}\}$ and $\{y'_{i,n}\}$ be the points in
$\cE_j$ satisfying (\ref{Large pair}) and (\ref{Iterated distance}). We will omit here the subscript $n$ and just write $y_{i,n}=y_i$ and $y'_{i,n}=y_i'$.
 In the case when these sets of preimages are finite, let $k$ be the cardinality of $\{y_i\}$ (so we have $\{y_1,...,y_k\}$ and $\{y'_1,...,y'_k\}$), and otherwise set $k=\infty$.
Then we can write
\begin{eqnarray*}
\big|\cL_z^{j,n}g(x)-\cL_z^{j,n}g(x')\big|\\=
\big|\sum_{t=1}^k\big(e^{S_{j,n}f(y_t)+zS_{j,n}u (y_t)}g(y_t)-
e^{S_{j,n}f(y_t')+zS_{j,n}u(y_t')}g(y_t')\big)\big|
\\
\leq\sum_{t=1}^ke^{S_{j,n}f(y_t)+\Re(z)S_{j,n}u(y_t)}
|e^{i\Im(z)S_{j,n}u(y_t)}g(y_t)-e^{i\Im(z)S_{j,n}u(y_t')}g(y_t')|\\+
\sum_{t=1}^k|e^{i\Im(z)S_{j,n}u(y_t')}g(y_t')|
|e^{S_{j,n}f(y_t)+\Re(z)S_{j,n}f(y_t)}-e^{S_{j,n}f(y_t')+\Re(z)S_{j,n}u(y_t')}|\\:=
I_1+I_2.
\end{eqnarray*}
In order to estimate $I_1$, observe that for any $1\leq t\leq k$,
\begin{eqnarray*}
|e^{i\Im(z)S_{j,n}u(y_t)}g(y_t)-e^{i\Im(z)S_{j,n}u(y_t')}g(y_t')|\\\leq
|g(y_t)|\cdot|e^{i\Im(z)S_n^\om u(y_t)}-e^{i\Im(z)S_{j,n}u(y_t')}|+
|g(y_t)-g(y_t')|:=J_1+J_2.
\end{eqnarray*}
By the mean value theorem and then by (\ref{Regularity of the functions}),
\[
J_1\leq2\|g\|_\infty|\Im(z)|Q \rho^\al(x,x'),
\]
while by (\ref{Iterated distance}),
\[
J_2\leq v_{\al,\xi}(g)\rho^\al(y_t,y'_t)\leq
 v_{\al,\xi}(g)(\gam^{[n/n_0]})^{-\al}\rho^\al(x,x)
\] 
and it follows that
\begin{equation*}
I_1\leq \cL_{\Re(z)}^{j,n}\textbf1(x)
\big(2\|g\|_\infty|\Im(z)|Q+v_{\al,\xi}(g)
(\gam^{-\al[n/n_0]}\big)\rho^\al(x,x').
\end{equation*}
Finally, observe that
\begin{equation}\label{BB.1}
\cL_{\Re(z)}^{j,n}\textbf1(x)\leq 
\cL_0^{j,n}\textbf1(x)e^{|\Re(z)|\|S_{j,n}u\|_\infty},
\end{equation}
and so
\begin{eqnarray*}
I_1\leq \|\cL_0^{j, n}\textbf1\|_\infty e^{|\Re(z)|\|S_{j,n} u\|_\infty}\\\times
\big(2\|g\|_\infty|\Im(z)|Q)+v_{\al,\xi}(g)(\gam^{[n/n_0]})^{-\al}\big)
\rho^\al(x,x').
\end{eqnarray*}
Next, we  estimate $I_2$.  First, by the mean value theorem
and (\ref{Regularity of the functions}),
\begin{eqnarray*}
|e^{S_{j,n} j(y_t)+\Re(z)S_{j,n} u(y_t)}-
e^{S_{j,n} f(y'_t)+\Re(z)S_{j,n} u(y'_t)}|\leq(1+|\Re(z)|)Q\\\times
\max\{e^{S_{j,n} f(y_t)+\Re(z)S_{j,n} u(y_t)},
e^{S_{j,n} f(y'_t)+\Re(z)S_{j,n} f(y'_t)}\}
\rho^\al(x,x')
\end{eqnarray*}
 and therefore
\begin{eqnarray*}
I_2\leq(1+|\Re(z)|)\|g\|_\infty(\textbf L_{\Re(z)}^{j,n}\textbf1(x)+\cL_{\Re(z)}^{j,n}\textbf1(x'))
Q\rho^\al(x,x')\\\leq
2(1+|\Re(z)|)\|g\|_\infty\|\cL_0^{j,n}\textbf1\|_\infty e^{|\Re(z)|\|S_{j,n} u\|_\infty}     
Q\rho^\al(x,x')                                       
\end{eqnarray*}
where in the last inequality we used (\ref{BB.1}), yielding the first
statement of Proposition \ref{LY covering}.
Finally, by (\ref{BB.1}) we have
\[
\|\cL_z^{j,n}g\|_\infty\leq\|g\|_\infty\|\cL_0^{j,n}\textbf1\|_\infty\cdot e^{|\Re(z)|
\|S_{j,n}u\|_\infty}
\]
and the lemma follows from the above estimates, taking into account that 
\[
\|\cL_0^{j,n}\textbf{1}\|_\infty<\infty
\]
by our assumption.\qed

Next, for any $s>1$, consider the  real cones $\cC_{j,\bbR}\subset\cH_j$ given by 
\[
\cC_{j,\bbR}=\{g\in\cH_j: g\geq0\,\text{ and }g(x)\leq e^{sQd^\al(x,x')}\,\text{ if }d(x,x')<\xi\}
\]
where $Q$ is so that $f_j,u_j\in H_{Q,\al}$ and $s>1$ is a parameter which will be chosen later.
Let $\cC_j$ the canonical complexification of  $\cC_{j,\bbR}$ which is given by
\[
\cC_j=\{g\in\cH_j:\,\mu(g)\overline{\nu(g)},\,\,\,\forall\,\mu,\nu\in\cC_{j,\bbR}^*\}
\]
where $\cC_{j,\bbR}^*=\{\mu\in\cH_j^*:\,\mu(g)>0,\,\,\forall g\in\cC_{j,\bbR},\,g\not=0\}$ is the dual of the real cone $\cC_{j,\bbR}$. Recall that
the complex cone $\cC_j$ can also be written as
\[
\cC_j=\bbC'(\cC_j+i\cC_j)=\bbC'\{x+iy:\,x\pm y\in\cC_{j,\bbR}\}
\]
where $\bbC'=\bbC\setminus\{0\}$. 
 In the case when $\xi>1$, let $l_j$ be any probability measure on $\cE_j$, while in the case when $\xi\leq1$, set $l_j=\nu_j$, where $\nu_j$ is the measure from Assumption \ref{RPF 0 ASS}. Our main result here is the following

\begin{theorem}\label{Covering maps cones thm}
Suppose that Assumption \ref{CovAss} holds true with $\xi>1$, or that it holds true with $\xi\leq1$ and that Assumption \ref{RPF 0 ASS} holds true. Then
there exist positive constants $s_0,K,M,k_0, d_1$ and $r$ so that for any $s\geq s_0$ and $j\in\bbZ$:

(i) The cone $\cC_j$ is linearly convex and cones  $\cC_J$ and $\cC_j^*$ have bounded aperture:
 for any $g\in\cC_j$ and $\mu\in\cC^*$, 
\begin{equation}\label{aperture2.1}
\|g\|\leq K|l_j(g)|
\end{equation}
and
\begin{equation}\label{aperture2.2}
\|\mu\|\leq M|\mu(\textbf{1})|.
\end{equation}
In particular $l_j\in \cC_{j,\bbR}^*$.

(ii) The cone $\cC_j$ is reproducing:  for any $g\in\cH_j$ there exists $R(g)\in\bbC$ so that $|R(g)|\leq K_1\|g\|$
and 
\[
g+R(g)\in\cC.
\]

(iii) For any complex number $z\in B(0,r)$ and $k_0\leq k\leq 2k_0$ we have 
\[
\cL_z^{j,k}\bbC_j'\subset\bbC_{j+k}'
\]
and 
\[
\sup_{f,g\in\cC_j'}\del_{\cC_{j+k}}(\cL_z^{j,k}f,\cL_z^{j,k}g)\leq d_1
\]
where $\cC_j'=\cC_j\setminus\{0\}$.
\end{theorem}
Theoerm \ref{RPF} follows from Theorem \ref{YT cones thm}, exactly as in Section 3 of \cite{SeqRPF} and Chapter 5 of \cite{book}, and our main task in this Section is to prove Theorem 
\ref{YT cones thm}.

\begin{proof}
First, the proof of (\ref{aperture2.2}) proceeds exactly as the proof of Lemma 5.5.5 in \cite{book},
and the proof of Theorem \ref{Covering maps cones thm} (ii) proceeds exactly as the proof of Lemma 5.5.4 in \cite{book}.
Now we will prove (\ref{aperture2.1}).
Let $g\in\cC_{j,\bbR}$ and let $x,x'\in\cE_j$ be so that $d_j(x,x')<\xi$ and $g(x)\geq g(x')$.
The wither $g(x)=g(x')=0$ or $\min\big(g(x),g(x')\big)>0$ and then
\[
|g(x)-g(x')|=g(x)-g(x')\leq \big(e^{sQd^\al(x,x')}-1\big)g(x')
\]
and so by the mean value theorem,
\[
|g(x)-g(x')|\leq sQ e^{sQ\rho^\al(x,x')}\rho^\al(x,x')\|g\|_\infty.
\]
Reversing the roles of $x$ and $x'$, we conclude that this inequality holds true for any 
$x$ and $x'$ such that $\rho(x,x')<\xi$, and therefore 
\begin{equation}\label{v est}
v_{\al,\xi}(g)\leq sQ e^{sQ\xi^\al}\|g\|_\infty\leq sQe^{sQ}.
\end{equation}

Next, in the case when $\xi>1$ we clearly have
\[
\|g\|_\infty\leq e^{sQ}l_j(g)
\]
and so by (\ref{v est}) and Lemma \ref{Lemma 5.3}  we can take $K=2\sqrt2 (1+sQe^{sQ})e^{sQ}$.
In the case when $\xi\leq1$ we proceed in a different way, relying on Assumption \ref{RPF 0 ASS}. We first claim that there exists a constant $c>0$ so that for any $j$ and $x_0\in\cE_j$ we have
\begin{equation}\label{Lower}
\nu_j(B_j(x_0,\xi))\geq c
\end{equation}
where $B_j(x_0,\xi)=\{x\in\cE_j:\,d_j(x,x_0)<\xi\}$.
Relying on (\ref{Lower}), for any $g\in\cC_{j,\bbR}$ we have
\[
\|g\|_\infty\leq\sup_{x}\frac{e^{sQ}}{\nu_j(B_j(x,\xi))}\int_{B_j(x,\xi)}gd\nu_j\leq
e^{sQ}c^{-1}\nu_j(g)
\]
and so by Lemma \ref{Lemma 5.3} we can take $K=2\sqrt2(1+sQe^{sQ})e^{2sQ}c^{-1}$.
Now we will prove (\ref{Lower}). The triplet $(\la_j,h_j,\nu_j)$ from Assumption \ref{RPF 0 ASS} satisfy 
\[
\la_{j,n}\nu_j=(\cL_0^{j,n})^*\nu_{j+n}
\]
and so with $B=B_j(x_0,\xi)$, 
\[
\nu_j(B(x_0,\xi))=\frac{\int\cL_0^{j,m_0}\bbI_B d\nu_{j+_m}}{\la_{j,m_0}}\geq 
C^{-m_0}\min_x\cL_0^{j,m_0}\bbI_B(x) 
\]
where $\bbI_B$ is the indicator function of $B$ and $C>0$ is an upper bound of the sequence $\{\la_j\}$.
Now, by (\ref{TopEx}) for any $x\in\cE_{j+m_0}$ there exists $y=y_x\in\cE_j$ so that $T_j^{m_0}y_x=x$ and hence 
\[
\min_x\cL_0^{j,m_0}\bbI_B(x) \geq e^{-m_0\|f\|_\infty}
\]
where $\|f\|_\infty=\sup_j\|f_j\|_\infty<\infty$. Now we can take $c=C^{-m_0}e^{-m_0\|f\|_\infty}$.

Now we will prove Theorem \ref{Covering maps cones thm} (iii). We first claim that for any $j$ and $n\geq0$ we have
\begin{equation}\label{1}
\cL_0^{j,n_0+n}\cC_{j,\bbR,s}\subset\cC_{j+n_0+n,\bbR,s'}
\end{equation}
where $s'=s\gam^{-1}+1$.
Indeed, let  $x,x'\in\cE_{j+n_0+n}$ be so that  $d(x,x')=d_{j+n_0+n}(x,x')<\xi$, 
and let $\{y_i\}$ and $\{y'_i\}$ be their preimages under $T_{j}^{n+n_0}$ so that for each $i$,
\[
d_j(y_i,y_i')\leq\gam^{-[\frac{n+n_0}{n+0}]}d_{j+n+n_0}(x,x')\leq \gam^{-1}d_{j+n+n_0}(x,x').
\]
Then, since $\{f_j\}\in H_{\al,Q}$, 
for any $g\in\cC_{j,\bbR,s}$ we have
\begin{eqnarray*}
\cL_0^{j,n_0+n}g(x)=\sum_{i}e^{S_{j,n_0+n}f(y_i)}g(y_i)\leq 
e^{Qd^\al(x,x')+sQ\gam^{-1}d^\al(x,x')}\sum_{i}e^{S_{j,n_0+n}f(y'_i)}g(y'_i)\\
=e^{(1+\gam^{-1}s)Qd^\al(x,x')}\cL_0^{j,n_0+n}g(x').
\end{eqnarray*}
Now we will fix some $s>\frac1{1-\gam^{-1}}$. We claim next that there exists constants $C>0$ and
$k_0\geq n_0$ so that for any $k\geq k_0$,  an integer $j$, a function $g\in\cC_{j,\bbR}$ and $x,x'\in\cE_{j+k}$ we have
\[
\cL_0^{j,k}g(x)\leq C \cL_0^{j,k}g(x').
\] 
In the case when $\xi>1$ this already follows from (\ref{1}), where we can just take $k_0=n_0$. In the case when $\xi\leq1$ using Assumption \ref{RPF 0 ASS} we have 
\[
\frac{\cL_0^{j,n}g(x)}{\cL_0^{j,n}g(x')}=\frac{\la_{j,n}^{-1}\cL_0^{j,n}g(x)}{\la_{j,n}^{-1}\cL_0^{j,n}g(x)}
\leq\frac{h_{j+n}(x)+C\del_n\ka(g)}{h_{j+n}(x')-C\del_n\ka(g)}
\]
where $\ka(g)=\frac{\|g\|}{\nu_j(g)}$, $C>0$ is some constant and
\[
\del_n=\sup_{j\in\bbZ}\left\|\frac{\cL^{j,n}_0}{\la_{j,n}}-\nu_j\otimes h_{j+n}\right\|_\al.
\]
Since $h_{j+n}(x)\in[a,b]$ for some constants $a,b>0$ which do not depend on $j$ and $n$, and $|\ka_j(g)|\leq K$ we obtain that there exist $k_0$ and $C$ which satisfy there required conditions. 

The last step of the proof is the following estimate, whose proof proceeds exactly as the proof of Lemma 5.8.1 in \cite{book}: for any nonzero $g\in\cC_{j,\bbR}$ and $\nu\in\cC_{j+m,\bbR}^*$  we have
\begin{equation}\label{2}
|\nu(\cL_z^{j,m})-\nu(\cL_0^{j,m})|\leq \nu(\cL_0^{j,m})|z|\big(2\|S_{j,m}u\|_\infty+3(s-1)^{-1}\big).
\end{equation}
Note that in \cite{book} we worked under the assumption that the maps $T_j$ have finite degree, but in the proof of
Lemma 5.8.1 we only used (5.17) and (5.18) from \cite{book}, and it is clear that we can use (\ref{Large pair}) and (\ref{Iterated distance}) instead. Note also that in \cite{book} we compared between $\cL_z^{j,m}$ and 
$\cL_{\Re(z)}^{j,m}$, but the proof from there proceeds exactly the same when comparing between $\cL_z^{j,m}$ and $\cL_{\Re(z)}^{j,m}$, with the exception that each appearance of $\Im(z)$ should be replaced by $|z|$. 

Relying on the above, the proof of Theorem \ref{Covering maps cones thm} proceeds exactly as the proof of equality (5.9.1) in \cite{book}: we first use that the real projective diameter of 
\[
\cC_{j,\bbR,s'}\cap\{g: |g(x)|\leq C|g(y)|,\,\,\forall\,x,y\}
\]
inside $\cC_{j,\bbR,s}$
does not exceed $d_0:=2\ln\big(\frac{s+s'}{s-s'}C\big)$ (see Lemma 5.7.3 in \cite{book}), and then use (\ref{2}) and Theorem A.2.4 from Appendix A of \cite{book} (which was used after equality (\ref{Comp1})) to derive that for some $k_0$, $r>0$ and all $j$, $z\in \bbC$ with $|z|<r$  and $k_0\leq k\leq 2k_0$ we have
\begin{equation*}
\cL_z^{j,k}\cC_j'\subset\cC_{j+k}'
\end{equation*}
and
\begin{equation*}
\sup_{f,g\in\cC_j'}(\cL_z^{j,k}f,\cL_z^{j,k}g)\leq d_0+6|\ln(1-\del)|.
\end{equation*}
where $\del\in(0,1)$ is some constant.
\end{proof}


\subsection{Examples}\label{examples}
We have listed several examples when Assumptions \ref{CovAss} (with $\xi\leq 1)$ and Assumption \ref{RPF 0 ASS} hold true,
and here we will present some models in which Assumption \ref{CovAss} holds true with $\xi>1$.
First, Assumption \ref{CovAss} holds true with $\xi>1$ and $n_0=1$ when for each $j$ we have $\cE_j=[0,1)$ and there exists a finite or countable collection of disjoint intervals $I_{j,n}=[a_{j,n},b_{j,n})$ whose union cover $\cE_j$, and $T_j|I_{j,n}:I_{j,n}$ is onto $[0,1)$ and satisfies that 
\[
|T_j(x)-T_j(y)|\geq\gam|x-y|
\] 
for any $x,y\in I_{j,n}$, where $\gam>1$ is a constant which is independent of $n$. These are exactly the type of maps considered in Section 3 of \cite{Conze}.

Consider now the Gauss map $T:\cE\to\cE$ which is given by $Tx=1/x-[1/x]$, where $\cE=[0,1]\setminus\bbQ$, and for each $j$ set $T_j=T$ and $\cE_j=\cE$.
Then Assumption \ref{CovAss} holds true with $n_0=2$ and $\xi>1$. Indeed, the preimages of any two points $x,y$ are paired as follows:
\[
T^{-1}\{x\}=\left\{x_n:\,n\in\bbN\right\}\,\,\text{ and }\,\,
T^{-1}\{y\}=\left\{y_n:\,n\in\bbN\right\}
\] 
where for any $z\in(0,1)$ and $n\in\bbN$ we set $z_n=\frac1{z+n}$.
When $n>1$ then 
\[
\left|\frac1{x+n}-\frac1{y+n}\right|\leq\frac{|x-y|}{n^2}\leq \frac14|x-y|.
\]
For $n=1$ we have 
\[
\min(x_1,y_1)>\frac12
\]
and 
\[
|x_1-y_1|<|x-y|.
\]
The preimages of $x_1$ and $y_1$ satisfy 
\[
\left|\frac1{x_1+n}-\frac1{y_1+n}\right|\leq (2/3)^2|x_1-y_1|\leq  (2/3)^2|x-y|
\] 
and so we can take $\gam=9/4=\min(4,9/4)$. 

The classical dynamical system related to continued fractions is the one generated by $T$ together with its unique absolutely continuous probability measure, whose density is given by $h(x)=\frac1{\ln 2}(1+x)^{-1}$. This corresponds to considering the transfer operator defined by the formula
\[
\cL_0g(x)=\sum_{n=1}^{\infty}\frac1{(x+n)^2}g(\frac1{x+n})=\sum_{n=1}^{\infty}e^{-2\ln x_n}g(x_n).
\]
In order to apply our results for this transfer operator, we need to verify that the function $f(y)=f_j(y)=-2\ln y$ satisfies (\ref{Regularity of the functions}). Set $\del=\gam^{-1}=\frac{4}9$. Denote by $x_{n_1,...,n_k}$ the inverse image by $T$ of $x_{n_1,...,n_{k-1}}$ which comes from the interval $(\frac1{n_k+1},\frac1{n_k})$, namely the  $x_{n_1,...,n_k}$'s are defined by the recursion formula $x_{n_1}=\frac1{x+n_1}$ and 
\[
x_{n_1,...,n_k}=\frac1{x_{n_1,...,n_{k-1}}+n_k}.
\]
Then for any $0\leq j<k$ we have $T^jx_{n_1,...,n_k}=x_{n_1,...,n_{k-j}}$.
The fact that the function $f(y)=-2\ln y$ satisfies (\ref{Regularity of the functions}) follows from the following claim: for $x,y\in\cE$, $k\geq0$ and $n_1,...,n_{2k+1},n_{2k+2}\in\bbN$ we have 
\[
\sum_{j=1}^{2(k+1)}\big|\ln(x_{n_1,...,n_j})-\ln(y_{n_1,...,n_j})\big|\leq 2(1+\del+\del^2+...+\del^{k})|x-y|.
\]
Indeed, the first two summands satisfy
\begin{eqnarray*}
|\ln(x_{n_1,n_2})-\ln(y_{n_1,n_2})|=|\ln(x_{n_1}+n_2)-\ln(y_{n_1}+n_2)|\leq\frac{|x_{n_1}-y_{n_1}|}{n_2}\leq |x-y|
\end{eqnarray*}
and 
\begin{eqnarray*}
|\ln(x_{n_1})-\ln(y_{n_1})|=|\ln(x+n_1)-\ln(y+n_1)|\leq\frac{|x-y|}{n_2}\leq |x-y|.
\end{eqnarray*}
Let $j>2$. Then
\begin{eqnarray*}
|\ln(x_{n_1,...,n_{2j+2}})-\ln(y_{n_1,...,n_{2j+2}})|=
|\ln(x_{n_1,...,n_{2j+1}}+n_{2j+2})\\-\ln(y_{n_1,...,n_{2j+1}}+n_{2j+2})|\leq \frac{|x_{n_1,...,n_{2j+1}}-y_{n_1,...,n_{2j+1}}|}{n_{2j+2}}
\end{eqnarray*}
and 
\begin{eqnarray*}
|\ln(x_{n_1,...,n_{2j+1}})-\ln(y_{n_1,...,n_{2j+1}})|=
|\ln(x_{n_1,...,n_{2j}}+n_{2j+1})\\-\ln(y_{n_1,...,n_{2j}}+n_{2j+1})|\leq \frac{|x_{n_1,...,n_{2j}}-y_{n_1,...,n_{2j}}|}{n_{2j+1}}.
\end{eqnarray*}
By (\ref{Iterated distance}) applied with $n_0=2$ we have
\[
\max\Big(|x_{n_1,...,n_{2j}}-y_{n_1,...,n_{2j}}|, |x_{n_1,...,n_{2j+1}}-y_{n_1,...,n_{2j+1}}|\Big)\leq\gam^j|x-y|
\]
which completes the proof of our claim.

Now we will describe the setup of non-stationary full-shifts, for which our results also hold true.
Let $d_j\in\bbN\cup\{\infty\}$ be a family of generalized numbers so that $d_j\geq2$, and set 
\[
C_j=\{k\in\bbN:\,k<d_j\}.
\]
For each $j$ set 
\[
\cE_j=\{(x_{j+m})_{m=0}^{\infty}\in\bbN^\bbN:\,x_{j+m}\in C_{j+m}\}
\]
and let $d_j(x,y)=2^{-\min\{n\geq0:\, x_{j+n}\not=y_{j+n}\}}$ be the dynamical distance on $\cE_j$.
For each integer $j$ define $T_j:\cE_j\to\cE_{j+1}$ by 
\[
T_j(x_{j},x_{j+1},x_{j+2},..)=(x_{j+1},x_{j+2},...).
\]
We also set $n_0=1$, $\gam=2$ and $\xi=2$. Then Assumption \ref{CovAss} holds true with these $\xi,n_0$ and $\gam$.


\begin{thebibliography}{Bow75}

\bibliographystyle{alpha}
\itemsep=\smallskipamount















\bibitem{Jon}
J. Aaronson, M. Denker, {\em Local Limit Theorems for Gibbs-Markov Maps}, Stoch. Dyn. 1 (2001), 193-237.


\bibitem{Aimino}
R. Aimino, M. Nicol and S. Vaienti. {\em Annealed and quenched limit theorems for random expanding dynamical systems}, Probab. Th. Rel. Fields 162, 233-274, (2015).


\bibitem{Arno}
P. Arnoux and A.Fisher, {\em Anosov families, renormalization and non-stationary
subshifts}, Erg. Th. Dyn. Syst., 25 (2005), 661-709.



\bibitem{Bir}
G. Birkhoff, {\em Extension of Jentzsch's theorem},
Trans. A.M.S. 85 (1957), 219-227.



 
 

 

\bibitem{Bus}
P.J. Bushell, {\em The Cayley-Hilbert metric and positive operators}, Linear Alg. and Appl. 84 (1986), 271-281.




\bibitem{Conze}
J-P Conze and A. Raugi, {\em Limit theorems for sequential expanding dynamical
systems}, AMS 2007.










\bibitem{Dub1}
L. Dubois, {\em Projective metrics and contraction principles for
complex cones}, J. London Math.
Soc. 79 (2009), 719-737.

\bibitem{Dub2}
L. Dubois, {\em An explicit Berry-Esseen bound for uniformly expanding maps on the interval}, Israel J. Math. 186 (2011), 221-250.



\bibitem{GH}
Y. Guivar\'ch and J. Hardy, {\em Th\'eor\`emes limites pour une classe de cha\^ines de Markov et
applications aux diff\'eomorphismes d'Anosov}, Ann. Inst. H. Poincar\'e Probab. Statist. 24 (1988), no. 1, 73-98.


\bibitem{book}
Y. Hafouta and Yu. Kifer, {\em Nonconventional limit theorems and random dynamics}, 
World Scientific, Singapore, 2018.






\bibitem{HafMD}
Y. Hafouta, {\em Nonconventional moderate deviations theorems and exponential concentration inequalities}, to appear in  Ann. Inst. H. Poincar\'e Probab. Statist.

\bibitem{Annealed}
Y. Hafouta, {\em Limit theorems for some skew products with mixing base maps},
preprint, arXiv:1808.00735.

\bibitem{HafEdge}
Y. Hafouta, {\em Asymptotic moments and Edgeworth expansions for some processes in random dynamical environment}, arXiv preprint 1812.06924.

\bibitem{SeqRPF}
Y. Hafouta, {\em A sequential RPF theorem and its applications to limit theorems for time dependent dynamical systems and inhomogeneous Markov chains}, arXiv preprint 1903.04018, 49 pages.

 
\bibitem{HH}
H. Hennion and L. Herv\'e, {\em Limit Theorems for Markov Chains and Stochastic Properties of Dynamical
Systems by Quasi-Compactness}, Lecture Notes in Mathematics vol. 1766, Springer, Berlin, 2001.


  








\bibitem{Kato}
T. Kato, {\em Perturbation theory for linear operators. Classics in Mathematics}. 
Springer-Verlag, Berlin, 1995.

\bibitem{Kifer-1996}
Yu. Kifer, {\em Perron-Frobenius theorem, large deviations, and random perturbations
in random environments},
 Math. Z. 222(4) (1996), 677-698.

\bibitem{Kifer-1998}
Yu. Kifer, {\em Limit theorems for random transformations and processes in random
environments},
Trans. Amer. Math. Soc. 350 (1998), 1481-1518.

\bibitem{Kifer-Thermo}
 Yu. Kifer, {\em Thermodynamic formalism for random transformations revisited},
Stoch. Dyn. 8 (2008), 77-102.





\bibitem{Liverani}
C. Liverani, {\em Decay of correlations}, Ann. Math. 142 (1995), 239-301.



\bibitem{MSU}
V. Mayer, B. Skorulski and M. Urba\'nski, {\em Distance expanding random mappings,
thermodynamical formalism, Gibbs measures and fractal geometry},
Lecture Notes in Mathematics, vol. 2036 (2011), Springer.
 




\bibitem{Nicol}
M. Nicol, A. Torok, S. Vaienti, {\em Central limit theorems for sequential and random intermittent dynamical systems}, 
Ergodic Theory Dynam. Systems, 38, pp. 1127-1153, 2016.




\bibitem{Rug}
H.H. Rugh, {\em Cones and gauges in complex spaces: Spectral gaps and complex
 Perron-Frobenius theory}, 
Ann. Math. 171 (2010), 1707-1752.



\bibitem{Sarig1}
O. Sarig, {\em Thermodynamic formalism for countable Markov shifts}, Ergodic Theory Dynam.
Systems 19, 1565-1595, 1999.

\bibitem{Sarig2}
O. Sarig, {\em Existence of Gibbs measures for countable Markov shifts}, Proc. Amer. Math.
Soc. 131 (2003), no. 6, 1751-1758.






\bibitem{Viv1}
V\'eronique Maume-Deschamps, {\em Correlation decay for Markov maps on a countable state space}, 
Erg. Th. Dyn. Syst. 21, 165-196 (2001)

\bibitem{Viv2}
V\'eronique Maume-Deschamps, {\em Projective metrics and mixing properties on towers}, Trans. Amer. Math. Soc. 353, 3371-3389 (2001)




\bibitem{Y1}
L.S. Young, {\em Statistical properties of dynamical systems with some hyperbolicity}, 
 Ann. Math. 7 (1998) 585-650.

\bibitem{Y2}
L.S. Young, {\em  Recurrence time and rate of mixing}, Israel J. Math. 110 (1999) 153-88.










\end{thebibliography}
\end{document}